\documentclass[11pt]{article}

\usepackage{amssymb}
\usepackage{amsmath}
\usepackage{bbm}
\usepackage{multirow}	
\usepackage{enumerate}

\usepackage{color}
\usepackage{hyperref}

\usepackage{xcolor}
\hypersetup{
    colorlinks,
    linkcolor={red!50!black},
    citecolor={blue!50!black},
    urlcolor={blue!80!black}
}

\usepackage{graphicx}

\usepackage{epstopdf}
\usepackage{centernot}

\numberwithin{equation}{section}

\usepackage{amsthm}
\newtheorem{theorem}{Theorem}[section]

\newtheorem{lemma}{Lemma}[section]

\newtheorem{definition}{Definition}[section]
\newtheorem{remark}{Remark}[section]

\newcommand{\norm}[1]{{\left\| \, #1 \, \right\|}}
\newcommand{\nnorm}[1]{{\left\vert\kern-0.25ex\left\vert\kern-0.25ex\left\vert \, #1 \, 
    \right\vert\kern-0.25ex\right\vert\kern-0.25ex\right\vert}}
\newcommand{\C}{\mathbb{C}}
\newcommand{\Z}{\mathbb{Z}}
\newcommand{\R}{\mathbb{R}}
\newcommand{\T}{\mathbb{T}}
\newcommand{\F}{\mathcal{F}}

\newcommand{\W}{\mathcal{W}}
\newcommand{\Tau}{\mathcal{T}}

\def\*#1{\mathbf{#1}}

\newcommand{\ip}[1]{{\left< \, #1 \, \right>}}

\newcommand{\Lt}{\mathrm{L}^2(\mathbb{T}^2)}
\newcommand{\Hs}[1]{\mathrm{H}^{#1}(\mathbb{T}^2)}\newcommand{\spec}[2]{\sigma_{\mathtt{#1}} \big( #2 \big) }
\newcommand{\tbs}{\overline{T}^*}
\newcommand{\lnu}{\lambda^{(\nu)}}
\newcommand{\unu}{u^{(\nu)}}

\usepackage{mathrsfs} 
\newcommand{\RR}{\mathscr{R}}
\renewcommand{\Re}{\operatorname{Re}}
\renewcommand{\Im}{\operatorname{Im}}
\renewcommand{\hat}{\widehat}

\def\thanksProject{This work was partially supported by the Natural Sciences and Engineering Research Council of Canada and Simon Fraser University.}
                                                   
\def\thanksJavierNilima{Department of Mathematics, Simon Fraser University, Burnaby, BC, V5A 1S6, Canada 
  ({\tt javiera@sfu.ca, nigam@math.sfu.ca}).}
             
\begin{document}

\title{Characterization of singular flows of zeroth-order pseudo-differential operators via elliptic eigenfunctions: a numerical study\thanks{\thanksProject}}

\author{{\sc Javier A. Almonacid}\thanks{\thanksJavierNilima} \quad 
        {\sc Nilima Nigam}\footnotemark[2]}
       
\date{}

\maketitle

\begin{abstract}
\noindent 
The propagation of internal gravity waves in stratified media, such as those found in ocean basins and lakes, leads to the development of geometrical patterns called ``attractors''. These structures accumulate much of the wave energy and make the fluid flow highly singular. In more analytical terms, the cause of this phenomenon has been attributed to the presence of a continuous spectrum in some nonlocal zeroth-order pseudo-differential operators. 
In this work, we analyze the generation of these attractors from a numerical analysis perspective. First, we propose a high-order pseudo-spectral method to solve the evolution problem (whose long-term behaviour is known to be not square-integrable). Then, we use similar tools to discretize the corresponding eigenvalue problem. Since the eigenvalues are embedded in a continuous spectrum, we compute them using viscous approximations. Finally, we explore the effect that the embedded eigenmodes have in the long-term evolution of the system.
\end{abstract}



\section{Introduction}

% - Introduce typical PDE
% - prior work in physics/math
% - time-domain PDE v/s viscous resonances

% The outline is not required, but we show an example here.
%The paper is organized as follows. Our main results are in
%\cref{sec:main}, our new algorithm is in \cref{sec:alg}, experimental
%results are in \cref{sec:experiments}, and the conclusions follow in
%\cref{sec:conclusions}.

The propagation of internal waves in a stratified medium and their interaction with the surrounding topography is a phenomenon that has been highly studied at different levels: theoretically \cite{lm08,Ogilvie2005}, numerically \cite{Brouzet2016,GSP2008}, and experimentally, both in two dimensions \cite{Davis2020,Hazewinkel2008,Maas1997} and three dimensions \cite{DrijfhoutMaas2007,Pillet2018}. In this case, the fluid flow develops geometrical patterns that make the velocity field highly singular. These singularities are usually known as ``attractors''.

In smoother scenarios, such as when a container of fluid is vibrated, its response can be described using the eigenmodes (eigenfunctions) of the system. These manifest as large-scale standing waves, which in many cases are visible to the naked-eye \cite{Maas1997}. Here, the eigenmodes are assumed to be smooth quantities that completely describe the system (and many times they form a basis of the space on which the solution resides), and therefore, the solution to the modelling PDE is expected to be smooth. The eigenfrequencies (eigenvalues) thus form a countable set. However, for the case of internal waves, Maas \cite{Maas2005} claims that ``attractors, rather than eigenmodes'' drive the response to the system. Moreover, the fact that these attractors are present for a wide range of frequencies (cf. \cite{Maas1997}) suggests the presence of a \textit{continuous spectrum} and, perhaps, a lack of eigenmodes.

For similar problems, such as the propagation of inertial waves in rotating fluids, the relationship between the spectral and dynamical properties of the underlying differential operator has already been established (see, e.g., \cite{Ralston1973} for the case of some specific two-dimensional containers). The mathematical characteristics of the generation and propagation of internal waves would not be discovered until recently. 

The development of attractors is a phenomenon that can now be explained by studying the flow of zeroth-order pseudo-differential operators \cite{cdv18,cs20,dz19}. In these works, different tools from the pseudo-differential calculus and microlocal analysis are used to obtain quantitative statements about the features of these irregular structures. Moreover, Colin de Verdi\`ere \& Saint-Raymond \cite{cdv18,cs20} were able to confirm many of the findings by Maas et al. \cite{Maas2005,Maas1997}. 

A common element in these cases is the study of the spatio-temporal equation
\begin{equation} \label{eq:main_pde}
i\partial_t u + P(x,D)u = f e^{-i\omega_0 t},
\end{equation}
where $i=\sqrt{-1}$, $f$ is a smooth right-hand side, $\omega_0 \geq 0$ is a forcing frequency, and $P(x,D)$ is a zeroth-order pseudo-differential operator that may arise from the manipulation of a standard set of fluid equations (we refer the reader to \cite{thesis} for more details on how to obtain equations of this type). Here, $P$ is assumed to be a bounded, self-adjoint operator satisfying some dynamical assumptions which would be lated relaxed (to a certain degree) by Colin de Verdi\`{e}re \cite{cdv18} and Dyatlov \& Zworski \cite{dz19} In particular, the study of \eqref{eq:main_pde} in \cite{cdv18, cs20, dz19} views the evolution as the flow on an energy manifold generated by the Hamiltonian vector field $H_p$ associated to the principal symbol of $P$, which is assumed to be Morse-Smale with no fixed points.

The main results in \cite{cdv18, cs20, dz19} are essentially the same but obtained with different methods. First, they confirm that internal-wave attractors form when $P$ has some continuous spectrum. More precisely, the spectrum of $P$ in a neighbourhood of $\omega_0$ is shown to be absolutely continuous with possibly finitely many embedded eigenvalue. The corresponding eigenmodes are analytic \cite{Wang2020}. Then, these works prove that the long-term evolution of the flow (where the attractors are fully developed) is not a square-integrable function, but a distribution living in Sobolev spaces $\rm{H}^s$ of negative order.

It is evident that the computation of the embedded eigenmodes of $P$ is a challenging problem.  Works such as the one by Rieutord et al. \cite{RGV01} suggest that a way to overcome this difficulty is to consider a regularized problem which introduces a small viscosity $\nu>0$. The interest is then to compute the eigenmodes of the viscous operator $P+i\nu\Delta$ (with $-\Delta$ the standard Laplacian) for very low viscosities. While this makes sense from a physical point of view, the mathematical picture is not too straightforward. 

On the one hand, the operator $P$ is of zeroth-order with a combination of continuous spectrum and embedded eigenvalues. On the other hand, $P+i\nu\Delta$ is a second-order operator with a purely discrete spectrum. Hence, the eigenvalues of $P+i\nu\Delta$ may not necessarily converge to those of $P$ as $\nu \to 0^+$. However, in some other contexts, there is mathematical evidence that this limit makes sense for some elements of the spectrum of $P$ (see, for instance, the viscous approximation of Pollicot-Ruelle resonances by Dyatlov \& Zworski \cite{dz15}). An answer to this problem (for the specific $P$ under study) was given by Galkowski \& Zworski \cite{gz19}. There, it is shown that the limit set of eigenvalues of $P+i\nu\Delta$ as $\nu \to 0^+$ gives a set of resonances that includes some of the embedded eigenvalues of $P$ in a neighbourhood of 0.

According to the above, we focus in this work on the development of numerical tools to approximate the solution to equations of the form \eqref{eq:main_pde} in the 2-torus and to compute eigenpairs of the operator $P+i\nu \Delta$. We will consider a particular class of zeroth-order operators, from which instances have been found in \cite{dz19}. Because of the absence of boundaries, different forms of $P$ will give rise to attractors with different shapes. In general, the distributional character of the evolution and the non-square-integrability for long times, as well as the embedding of the eigenvalues in the continuous spectrum will be the main challenges in this computational study.

Because the problems are posed on a periodic domain, we discretize the equations using a pseudo-spectral approach, which has proven to be a powerful tool to solve nonlocal problems in periodic domains (see, for instance, \cite{AntoineLorin2019,k2008,t00}). We then use these tools to explore the effect that the embedded eigenmodes have on the long-term evolution of the flow, for which we resort to the viscous approximation presented above. Some of these results will be explained using the microlocal analysis ideas from \cite{cs20,dz19}. In general, emphasis will be put into analyzing the results in frequency space, since this will give us insight into how smooth a function (or distribution) is (see, for instance, \cite{cgss2015,rnt12}).

The rest of this work is organized as follows. First, we finish this section with some notation that will be used throughout the paper. Next, in Section \ref{sec:background} we establish some analytic background to better understand the problem, as well as provide a simple example that allows some computations ``by hand''. Then, in Section \ref{sec:3}, we present the pseudo-spectral techniques and time stepping methods used to solve the evolution and eigenvalue problems of interest. In Section \ref{sec:singular_beh}, we study numerically the regularity the approximations and how they follow previous analytical findings. In addition, we explore how the attractors relate to the energy manifolds on which the flow takes place. Then, in Section \ref{sec:spectra} we analyze in detail the spectrum of $P+i\nu \Delta$ in the limit as $\nu \to 0^+$. In particular, we explore how the viscous approximations to the embedded eigenmodes of $P$ partially characterize the long-term evolution of the solution. Finally, in Section \ref{sec:conclusions}, we present some conclusions about this work.

\subsection*{Notation}

Throughout this paper, we denote by $\T^2 := \R^2 \backslash \Z^2$ the standard 2-torus, by $|\cdot|$ as either the modulus of a complex number or the Euclidean norm of a point in $\R^2$, and for any $\xi \in \R^2$, we define Japanese bracket $\ip{\cdot}$ as $\ip{\xi} := (1+|\xi|^2)^{1/2}$. 

The space of rapidly decaying functions on the integer lattice $\Z^2$ will be denoted by $\mathcal{S}(\Z^2)$, whereas $\mathcal{S}'(\Z^2)$ will denote the space of tempered distributions, that is, continuous linear functionals on $\mathcal{S}(\Z^2)$ (where the continuity is defined in term of the usual seminorm-induced topology). 
%We remark that elements $u \in \mathcal{S}'(\Z^2)$ are of the form
%\[
%\varphi \mapsto \ip{u,\varphi} := \sum_{\xi \in \Z^2} u(\xi) \varphi(\xi) \quad \forall \, \varphi \in \mathcal{S}(\Z^2),
%\]
%where $u:\Z^2 \to \C$ grows at most polynomially at infinity, i.e. there exists constants $M, C_{u,M} > 0$ such that
%\[
%|u(\xi)| \leq C_{u,M} \ip{\xi}^M \quad \forall \, \xi \in \Z^2.
%\]
Also, we define the space $C^{\infty}(\T^2) := \cap_{m \geq 1} C^m(\T^2)$, where $C^m(\T^2)$ is the space of $m$-times continuously differentiable periodic functions on $\T^2$. In turn, the space of periodic distributions on $\T^2$ will be denoted by $\mathcal{D}'(\T^2)$. Hence, we consider the (toroidal) Fourier transform $\F:C^\infty(\T^2) \to \mathcal{S}(\Z^2)$ defined as
\[
(\F f)(\xi) \equiv \hat{f}(\xi) := \dfrac{1}{2\pi} \int_{\T^2} f(x) \, e^{-ix\cdot \xi} \, dx.
\]
The operator $\F$ is a bijection \cite{rt10}. Its inverse $\F^{-1}:\mathcal{S}(\Z^2) \to C^\infty(\T^2)$ is given by
\[
(\F^{-1} g)(x) = \dfrac{1}{2\pi} \sum_{\xi \in \Z^2} g(\xi) \, e^{ix\cdot \xi}.
\]
Finally, for $s \in \R$, we shall consider the Sobolev space $\Hs{s}$ defined by
\begin{equation} \label{eq:def_Hs}
\Hs{s} := \left\{ f: \T^2 \to \R \hspace{0.5em} : \hspace{0.5em}  \norm{f}_s^2 := \sum_{\xi \in \Z^2} \ip{\xi}^{2s} \, |\hat{f}(\xi)|^2 < \infty \right\}.
\end{equation}
In particular, $\Hs{0}$ corresponds to the usual Lebesgue space $\mathrm{L}^2(\T^2)$.

%%%%%%%%%%%%%%%%%%%%%%%%%%%%%%%%%%%

\section{General background}
\label{sec:background}

We begin by describing some analytical concepts that will be used in this work, such as pseudo-differential operators and energy surfaces. We also describe the main objects of study in this work.

\subsection{Periodic pseudo-differential operators}
% this will be important! pseudos of order 0 need to be properly explained here.

Let $p(x,\xi)$ be a function with variables in the phase space $\T^2 \times \Z^2$. We say that $p(x,\xi) \in C^\infty(\T^2 \times \Z^2)$ if $p(\cdot,\xi) \in C^\infty(\T^2)$ for all $\xi \in \Z^2$.

\begin{definition}[Toroidal symbol class $S^m$] \label{def:symbol_space}
Let $m \in \R$. The toroidal symbol class $S^m(\T^2 \times \Z^2)$ (or simply $S^m$) consists of those functions $p(x,\xi) \in C^\infty(\T^2 \times \Z^2)$ satisfying
\[
\left| \Delta_\xi^\alpha \partial_x^\gamma p(x,\xi) \right| \leq C_{p,\alpha,\gamma,m} \ip{\xi}^{m-|\alpha|} \quad \forall \, x \in \T^2, \ \xi \in \Z^2,
\]
for some constant $C_{p,\alpha,\gamma,m} > 0$ and every multi-index $\alpha,\gamma$. In general, for a function $\sigma:\Z^n \to \C$ and a multi-index $\alpha$, the partial difference operator $\Delta_\xi^\alpha$ is defined as
\[
\Delta_{\xi_j}\sigma := \sigma(\xi+\delta_j) - \sigma(\xi), \quad \Delta_\xi^\alpha := \Delta_{\xi_1}^{\alpha_1} \Delta_{\xi_2}^{\alpha_2},
\]
where for $1\leq i,j\leq 2$, the quantity $\delta_j \in \mathbb{N}_0^2$ is given by a Kronecker delta, i.e. $(\delta_j)_i := \delta_{ij}$.
\end{definition}

\begin{definition}
Let $p(x,\xi) \in S^m$. For a function $v \in C^\infty(\T^2)$, we define the operator $P(x,D)$ by
\begin{equation} \label{eq:def_P_general}
P(x,D) v(x) := \dfrac{1}{2\pi}  \sum_{\xi \in \Z^2} e^{i x \cdot \xi} \, p(x,\xi) \hat{v}(\xi).
\end{equation}
We call $P(x,D)$ an $m$-th order pseudo-differential operator, and we call $p(x,\xi)$ the symbol of $P(x,D)$.
\end{definition}

Pseudo-differential operators appear as generalizations of linear differential operators. In particular, $D$ is just a notation for the gradient
operator. For a more in depth description of this notation, we refer the reader to \cite{Hormander}. From this definition, we immediately notice that if the symbol is independent of $x$, then
\begin{equation}
P(D)v(x) = \F^{-1} \Big[ p(\xi) \, \hat{v}(\xi) \Big](x).
\end{equation}

 We also need the concept of the {\it principal symbol} of $P(x,D)$ which is defined via the symbol $p(x,\xi)$. For a precise definition we refer the reader to \cite{zworskibook}. For our purpose: consider the symbol 
 $$a(x,\xi) \sim \sum_{j=0}^\infty |\xi|^{m-j}a_j(x,\xi), \ |\xi|>1, \ a_j(x,t\xi) = ta_j(x,\xi), \ t>0.$$ 
 The principal symbol of $a(x,D)$ is then $|\xi|^m a_0(x,\xi)$; we may think of this as the `leading order derivative term'.

The operator $P(x,D)$ is continuous and maps the space $C^\infty(\T^2)$ into itself. Furthermore, if there exists $C>0$ such that
\[
\left| \partial_x^\gamma p(x,\xi) \right| \leq C \quad \forall \, x \in \T^2, \ \xi \in \Z^2,
\]
and for all $|\gamma| \leq 3$, then the operator $P(x,D)$ extends to a bounded operator on $\mathrm{L}^2(\T^2)$ (cf. \cite[Theorem 9.1]{rt10}).

\subsection{Problem statement}

Let $D = (D_{x_1},D_{x_2})$, where $D_{x_j} = -i\partial_{x_j}$; $r>0$ and $\beta \in C^\infty(\T^2)$ be a purely real function. In what follows, $P(x,D):\Lt \to \Lt$ will denote the zeroth-order pseudo-differential operator given by
\begin{equation} \label{eq:def_P}
P(x,D) v(x) = \ip{D}^{-1}D_{x_2} v(x) - r\beta(x) v(x), \quad v \in \Lt.
\end{equation}
Using the representation \eqref{eq:def_P_general}, this operator can also be written as
\begin{equation} \label{eq:P_kn_rep}
P(x,D)v(x) = \dfrac{1}{2\pi} \sum_{\xi \in \Z^2} e^{i x\cdot\xi} \, \big[ \ip{\xi}^{-1}\xi_2 - r\beta(x) \big] \, \widehat{v}(\xi) .
\end{equation}

We will be interested in the discretization of two different (but related) problems. First, we consider the problem of finding a complex function $u = u(x,t)$, $x=(x_1,x_2) \in \T^2$, $t \geq 0$ such that
\begin{equation} \label{eq:evol_problem}
iu_t - Pu = f e^{-i\omega_0 t} \quad \text{in }\T^2 \times (0,\infty), \qquad u\big|_{t=0} = 0.
\end{equation}
where $f \in C^\infty(\T^2)$ and $\omega_0 \geq 0$ is a forcing frequency. Next, let $\Delta$ be the standard Laplacian operator and $\nu >0$ (parameter wich we shall call ``viscosity''). Then, we wish to find pairs $(\lambda^{(\nu)},u^{(\nu)})$ such that
\begin{equation} \label{eq:ev_problem}
(P(x,D) - \omega_0 + i\nu\Delta) u^{(\nu)}(x) = \lambda^{(\nu)} u^{(\nu)}(x), \quad x \in \T^2.
\end{equation} 

\subsection{Dynamical assumptions and the concept of energy manifolds} \label{sec:assumptions}

Let $\tbs \T^2$ denote the fibre-radially compactified cotangent bundle of $\T^2$. Following to \cite{cdv18,cs20,dz19}, we assume that $P(x,D)$ is a zeroth-order, self-adjoint operator (with respect to the usual $L^2$ inner product) with principal symbol $\bar{p} \in S^0(\tbs \T^2 \backslash \{0\})$ (the definition of this space is similar to that in Definition \ref{def:symbol_space}, but we refer to \cite{dzbook} for more details). % which is {\color{red} homogeneous of degree 0}.
It can be seen that for the particular choice of $P$ in \eqref{eq:def_P}, these assumptions are indeed satisfied (cf. \cite{dz19}).

The work by Dyatlov \& Zworski \cite{dz19} provides us with an statement that will be of utmost importance to understand the dynamics of the problem \eqref{eq:evol_problem}. Indeed, consider the quotient map for the $\R^+$ action $(x,\xi) \mapsto (x,s\xi)$, $s>0$,
\[
\kappa: \tbs \T^2 \to \partial \tbs \T^2.
\]
Then, the rescaled Hamiltonian vector field $|\xi|H_{\bar{p}}$, where
\begin{equation} \label{eq:def_hamiltonian}
H_{\bar{p}} := \sum_{j=1}^2 \partial_{\xi_j} \bar{p} \partial_{x_j} - \partial_{x_j} \bar{p} \partial_{\xi_j},
\end{equation}
commutes with the $\R^+$ action and the flow of
\begin{equation}
X := \kappa_*(|\xi|H_{\bar{p}}) \quad \text{is tangent to}\quad \Sigma := \kappa \big( \bar{p}^{\, -1}(\{0\}) \big).
\end{equation}
This crucial last statement will be addressed in more detail in Section \ref{sec:singular_beh}. We call $\Sigma$ the \textit{energy manifold}. To ensure that attractors are generated, one can make assumptions about the flow on $\Sigma$, but their verification can become rather difficult. Instead, we choose to study the energy manifold $\Sigma$ in itself. Let us consider a simple case where  $\Sigma$ can be characterized more explicitly. This will be useful to describe the exact location where the attractors are expected to appear.

\subsection{A simple example} \label{sec:example_cos}

Consider for $r>0$ the symbol and its corresponding principal symbol:
\begin{equation} \label{eq:p_example}
p(x,\xi) := \xi_2 \ip{\xi}^{-1} - r \cos(x_1), \quad \bar{p}(x,\xi) = \xi_2|\xi_2|^{-1} - r\cos(x_1).
\end{equation}
Then, according to \eqref{eq:def_hamiltonian}, the dynamical (Hamiltonian) equations are
\begin{equation} \label{eq:sys1}
\left\{ \begin{aligned}
&\frac{dx_1}{dt} = -r|\xi| \sin(x_1), \quad &&\frac{dx_2}{dt} = 0, \\
&\frac{d\xi_1}{dt} = \dfrac{\xi_1 \xi_2}{|\xi|^2}, \quad &&\frac{d\xi_2}{dt} = -\dfrac{\xi_1^2}{|\xi|^2}.
\end{aligned} \right.
\end{equation}
Given that the energy surface $\Sigma_0$ takes the form
\begin{equation} \label{eq:def_s0}
\Sigma_0 = \kappa\left( \left\{ (x_1,x_2,\xi_1,\xi_2) \in\tbs \T^2 \backslash \{0\}  : \dfrac{\xi_2}{|\xi|} = r\cos(x_1) \right\} \right),
\end{equation}
the system \eqref{eq:sys1} restricted to $\Sigma_0$ yields
\begin{equation} \label{eq:sys2}
\left\{ \begin{aligned}
&\frac{dx_1}{dt} = -r|\xi| \sin(x_1), \quad &&\frac{dx_2}{dt} = 0, \\
&\frac{d\xi_1}{dt} = r \dfrac{\xi_1}{|\xi|}\cos(x_1), \quad &&\frac{d\xi_2}{dt} = -\dfrac{\xi_1^2}{|\xi|^2}.
\end{aligned} \right.
\end{equation}
Here, we see that $\frac{d\xi_2}{dt}<0$, so $\xi_2$ must be decreasing as $t \to +\infty$. Given that $|\xi|$ is constant (thanks to the fibre-compactification), and assuming that $|\xi_1|>|\xi_2|$, we must have that $\xi_2 \to 0$. This, in turn, implies that $r \cos(x_1) \to 0$ (using \eqref{eq:def_s0}). Since $r>0$, this means that $x_1 \to \frac{\pi}{2}$ or $x_1 \to -\frac{\pi}{2}$.

\begin{figure}
\centering
\includegraphics[width=0.99\textwidth]{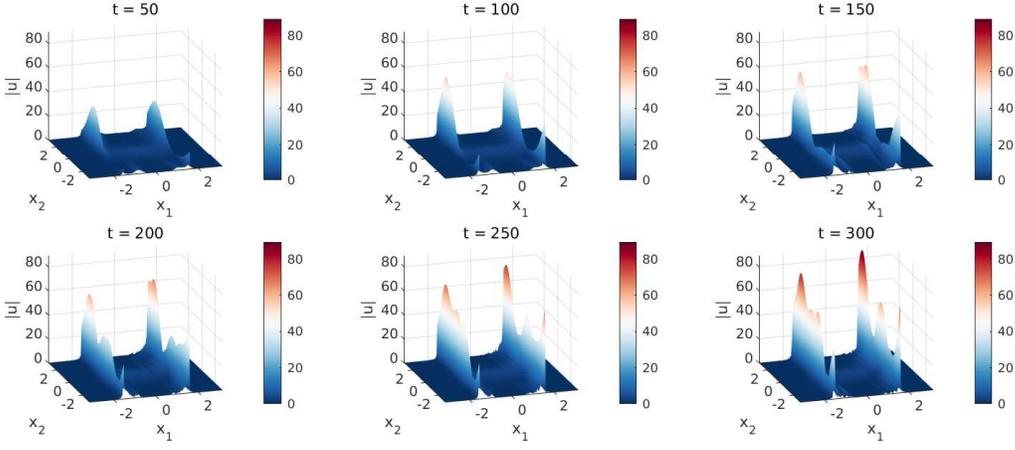}
\vspace*{-1em}
\caption{Development of attractors for the pseudo-differential operator given in \eqref{eq:p_example}. The solution was computed using the numerical techniques to be described in Section \ref{sec:3}. \label{fig:attractors}}
\end{figure}

To complete the description of the set of attracting Lagrangians $\Lambda_0^+$, we have to analyze the behaviour of $\xi_1$. Notice that, as $\xi_2 \to 0$, $|\xi_1|$ increases,  and therefore, if $\xi_1 < 0$, we will need $\frac{d\xi_1}{dt} <0$, and if $\xi_1 > 0$ we will need $\frac{d\xi_1}{dt} > 0$. First, linearizing \eqref{eq:sys2} around $x_1 = \pi/2$, $x_2 \in \mathbb{S}^1$, $\xi_2 = 0$, we notice that $\xi_1$ must be negative in order to get a sink of the dynamical system. Similarly, linearizing around $x_1 = -\pi/2$, $x_2 \in \mathbb{S}^1$ and $\xi_2=0$, we have that $\xi_1$ must be positive to get a source of the system. Therefore,
\begin{equation} \label{eq:attracting_lag}
\Lambda_0^+ = \left\{x_1 = \frac{\pi}{2}, \ x_2 \in \mathbb{S}^1, \ \xi_1 < 0, \ \xi_2 = 0 \right\} \cup \left\{x_1 = -\frac{\pi}{2}, \ x_2 \in \mathbb{S}^1, \ \xi_1 > 0, \ \xi_2 = 0 \right\}.
\end{equation}
This set of attractors is portrayed in Figure \ref{fig:attractors}. This was computed using the numerics to be described in the next section.

\section{Spectral discretization and convergence studies} \label{sec:3}

In both the evolution problem \eqref{eq:evol_problem} and the eigenvalue problem \eqref{eq:ev_problem}, the operators $P(x,D)$ and $P(x,D)+i\nu\Delta$, respectively, will be discretized using a pseudo-spectral approach. First, we discretize $\T^2$ using a mesh $\Tau_N$ containing $N$ points per direction (with $N$ even), that is,
\begin{equation} \label{eq:mesh}
\mathcal{T}_N := \left\{ x_j = (x_{j_1},x_{j_2}) \in \R^2 \ : \ x_{j_\iota} = 2\pi j_\iota/ N, \ j_\iota = -N/2,\dots,N/2-1 \right\}.
\end{equation}
Similarly, in frequency space, we consider a set $\W_N$ with $N$ wave-numbers per direction, that is,
\begin{equation} \label{eq:def_wavenumbers}
\W_N := \big\{k = (k_1,k_2) \in \Z^2 \ : \ k_\iota = -N/2, \dots, N/2-1 \big\}.
\end{equation}
Then, we approximate $u(x,t) \approx u_N(x,t)$, where
\[
u_N(x_j,t) = \dfrac{1}{N^2} \sum_{k_1=-N/2}^{N/2-1} \sum_{k_2=-N/2}^{N/2-1} \widehat{u}_N(k,t) e^{2\pi i (k\cdot x_j)/N}, \quad x_j \in \Tau_N, \ t > 0.
\]
and $\widehat{u}_N$ are Fourier coefficients computed using the discrete Fourier transform (DFT) defined as
\begin{equation} \label{eq:dft}
(F u_N)(k) \equiv \widehat{u}_N(k) := \sum_{j_1=-N/2}^{N/2-1}\sum_{j_2=-N/2}^{N/2-1} u_N(x_j) \, e^{-2\pi i (k\cdot j)/N}, \quad k \in \W_N.
\end{equation} 

It turns out that the representation \eqref{eq:P_kn_rep} is fundamental to understand how to transform the action of $P$ in frequency space. Indeed, the semi-discrete version of \eqref{eq:evol_problem} and the discrete version of \eqref{eq:ev_problem} become respectively:
\begin{equation} \label{eq:evol_fourier}
i \, \partial_t \widehat{u}_N - \ip{k}^{-1}k_2 \, \widehat{u}_N + r F \left( \beta \, F^{-1} \widehat{u}_N \right) = \widehat{f} e^{-i\omega_0 t}, \quad k \in \W_N, \ t>0,
\end{equation} 
and
\begin{equation} \label{eq:ev_fourier}
\left( \ip{k}^{-1}k_2 - \omega_0 - i\nu|k|^2 \right) \hat{u}_N^{(\nu)} - r F \left( \beta \, F^{-1} \hat{u}_N^{(\nu)} \right) = \lnu \unu_N, \quad k \in \W_N.
\end{equation}

To compute the error in the approximation, we consider a discrete version of the $\Hs{s}$-norm in \eqref{eq:def_Hs}. For a discrete function $u_N$ defined on $\Tau_N$, and $N\geq 2$ an even integer, we define:
\begin{equation} \label{eq:dnorm}
\nnorm{u_N}_s^2 := \dfrac{h^2}{N^2} \sum_{k_1=-N/2}^{N/2-1}\sum_{k_2=-N/2}^{N/2-1} \ip{k}^{2s} \, |\widehat{u}_N(k)|^2,
\end{equation}
where $h := 2\pi/N$ is the grid spacing. The scaling factor $h^2/N^2$ turns this discrete norm into a true approximation of $\norm{\cdot}_s$ as $N \to \infty$.

\begin{remark}
  Note that de-aliasing techniques have not being considered so far. On the one hand, it is known that pseudo-spectral methods amplify aliasing errors when the solution lacks regularity. On the other hand, standard tools such as the two-thirds rule or spectral viscosity methods (cf. \cite{BardosTadmor2015}) may not be the best alternative for this problem, as we are expecting singular solutions that can be no more regular than $\mathrm{L}^2(\T^2)$. The high frequencies will be key to describe these singularities. Therefore, removing them without careful thought might lead to an unrealistic smoothing of the dynamics. Thus, it still remains to determine a proper de-aliasing strategy for this problem.
\end{remark}

\subsection{The evolution problem}

We are interested in pairing the pseudo-spectral method in \eqref{eq:evol_fourier} with a high-order, one-step time discretization. The first option is to consider an exponential time-differencing fourth order Runge-Kutta method \cite{kt2005}, due to its demonstrated reliability in nonlinear problems \cite{k2008}. We notice in addition that, because the second term in \eqref{eq:evol_fourier} results in a diagonal matrix, its eigenvalues are precisely the diagonal entries, all of which satisfy $|\ip{\xi}^{-1}\xi_2| \leq 1$. This suggests that a simple fourth order Runge-Kutta method would also work for this problem, yielding a rather mild time step restriction: $\Delta t \leq 2.7$ (cf., e.g., \cite{t00}). 

\begin{figure}
\centering
\includegraphics[width=0.445\textwidth]{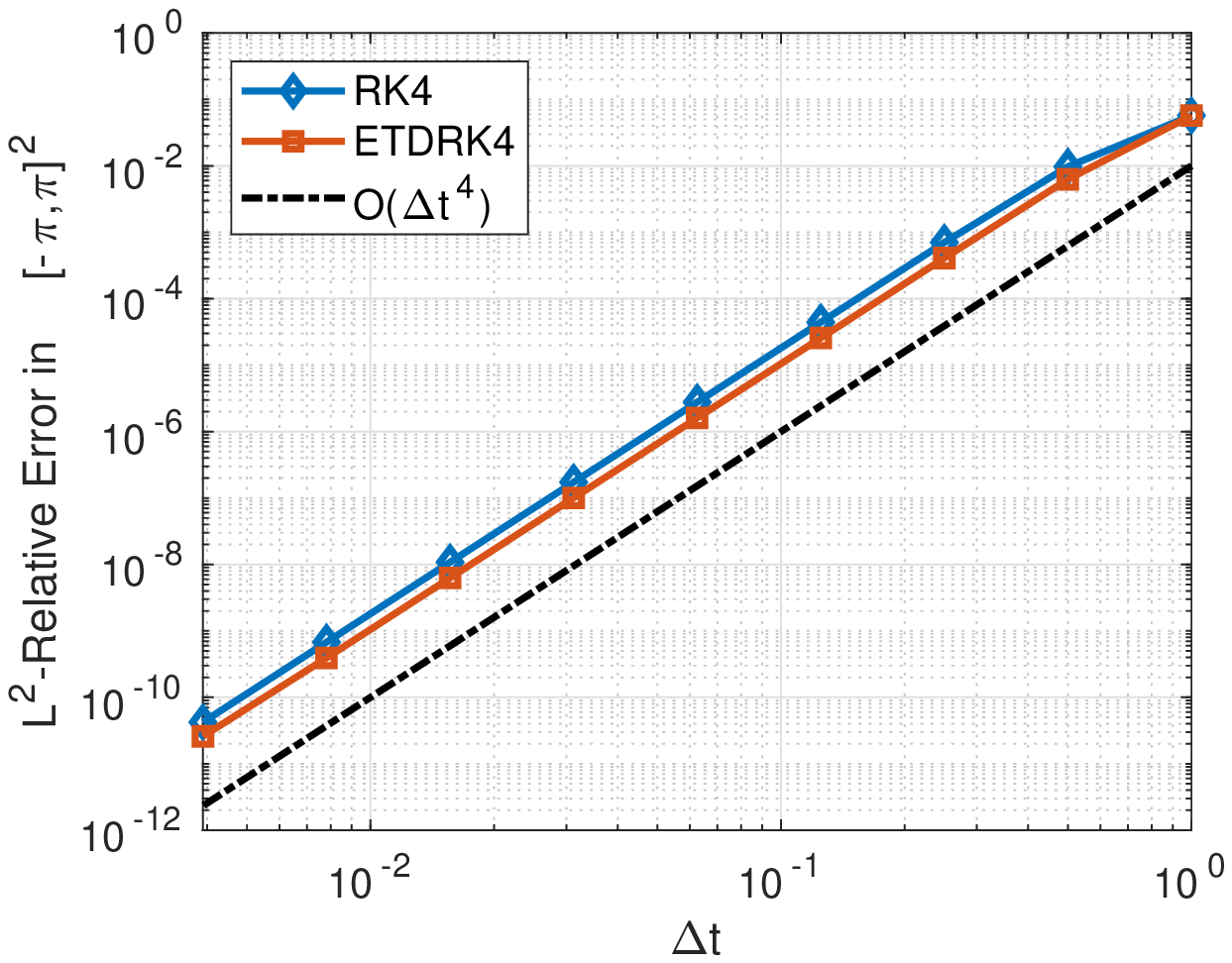}\hspace{12pt}
\includegraphics[width=0.445\textwidth]{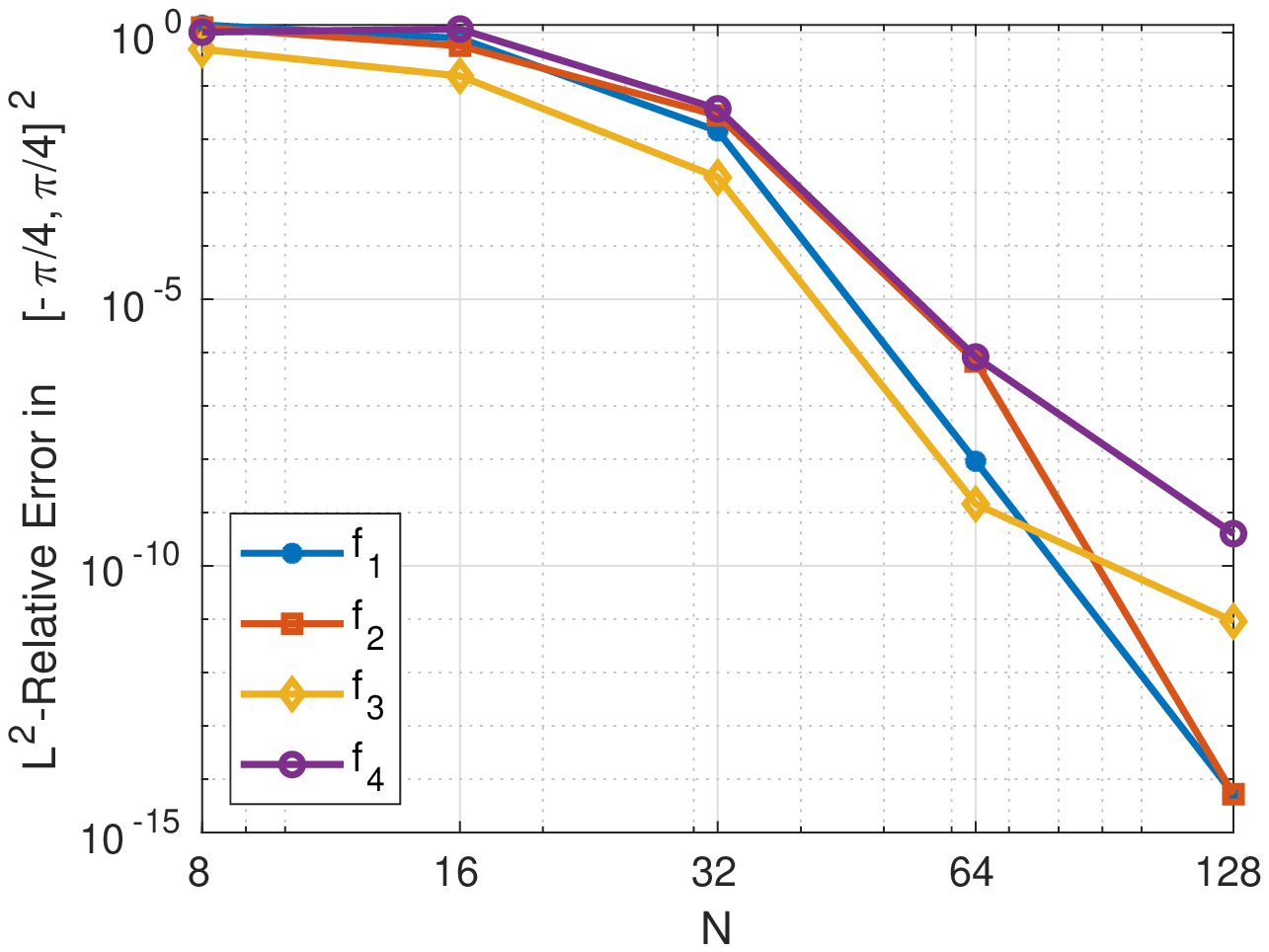}
\caption{Convergence studies. Left: fourth-order convergence in time. Right: local accuracy in space.\label{fig:convergence}}
\end{figure}

A first glimpse on how the numerics work was already shown in Figure \ref{fig:attractors}. Here we constructed the solution using $\omega_0 = 0$, $r=2$, $\beta(x) = \cos(x_1)$ and $f(x) = -5 \left( -|x|^2 + i(2x_1+x_2) \right)$. To study the convergence in time, we take the same data as before, but with $\omega_0 = 0.1$ and
\begin{equation} \label{eq:f4}
f(x) = -5\exp \Big( -3\big( (x_1+0.9)^2+(x_2+0.8)^2 \big) + i \big( 2x_1+x_2 \big) \Big).
\end{equation}
We fix a spatial mesh with $N = 64$ grid points per direction, and take several choices of time steps $\Delta t = 1, \frac{1}{2^1}, \frac{1}{2^2},\dots, \frac{1}{2^8}$. To test the robustness of our code, the peak of the Gaussian in  \eqref{eq:f4} has been purposely set closer to the attractor $x_1=-\pi/2$. We observe in Figure \ref{fig:convergence} (left) the expected fourth-order convergence, with the ETDRK4 method performing slightly better than the traditional RK4. Here, the error has been measured at the final time using the $\nnorm{\cdot}_0$ norm (over the whole domain) defined in \eqref{eq:dnorm}, and with respect to a more refined solution (same $N$ but with $\Delta t = 2^{-10}\cdot10^{-2} \approx 10^{-5}$).

The study of convergence in space requires a bit more thought. We notice in Figure \ref{fig:attractors} that the attractors at $x=\pm \pi/2$ manifest as singularities in the solution. Moreover, for fixed time $t>0$, there is no guarantee that $u(\cdot,t)$ is smoother than $\Lt$ (we go over this in more detail in Theorem \ref{th:uinf}). This tells us that global spectral accuracy cannot be expected. However, by looking again at Figure \ref{fig:attractors} we see that, in regions that are away from the attractors, the solution does appear to be smooth. Hence, we measure the $\rm L^2$-error with respect to a more refined solution ($N = 2^{10}$) in the domain $[-\pi/4,\pi/4]^2$. For this experiment we have fixed the time step at $\Delta t = 0.1$ and considered the ETDRK4 time-stepping method. Similarly to the study in time, we have taken $\omega_0 = 0.5$, $T=10$, $r=2$, $\beta(x) = \cos(x_1)$. We take values of $N = 2^3,2^4,\dots,2^7$ per spatial direction and also consider several source terms (with different levels of easiness of resolution):
\begin{align}
\label{eq:f1}f_1(x) &= \sin(x_1)\cos(2x_2), \\
f_2(x) &= \sin(x_1)\cos(2x_2) + \sin(5x_1)\cos(2x_2) + i\sin(5x_1)\cos(4x_2), \\
\label{eq:f3}f_3(x) &= \frac{1}{2}\exp\big( -2|x|^2 \big),
\end{align}
and $f_4$ as defined in \eqref{eq:f4}.  We show these results in Figure \ref{fig:convergence} (right), where we see rapid decay of the error in the region $[-\pi/4,\pi/4]^2$.

\subsection{The eigenvalue problem} \label{sec:3_ev_problem}

We solve the eigenvalue problem \eqref{eq:ev_fourier} using the matrix representation of the DFT and Matlab's built-in tool \texttt{eigs}, which uses an Arnoldi-type iteration to find a subset of eigenvalues of interest. In particular, we ask to compute the first $m$ eigenvalues closest to $\omega_0$ in the complex plane. 

\begin{figure}
\centering
\includegraphics[width=0.44\textwidth]{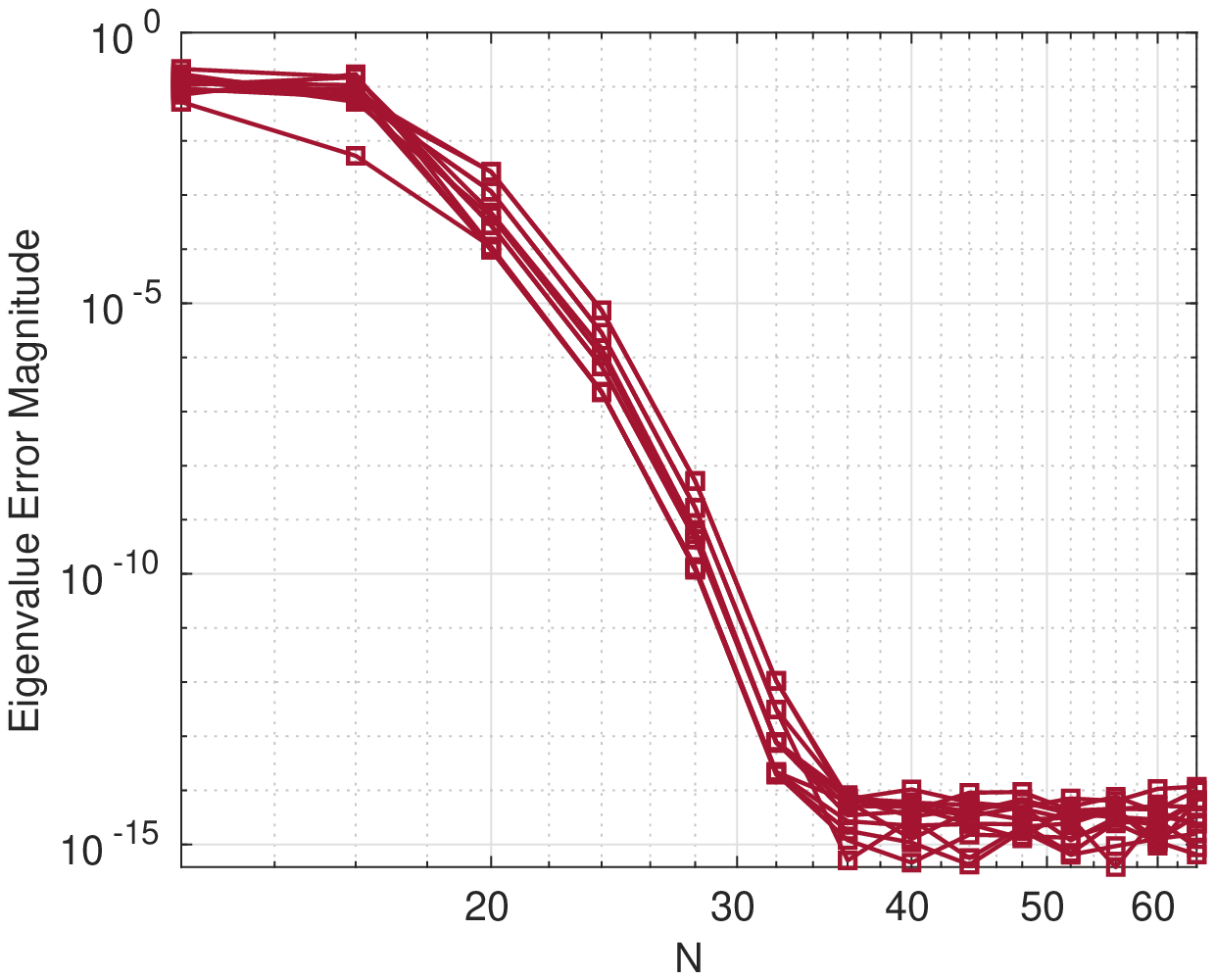} \ \ 
\includegraphics[width=0.44\textwidth]{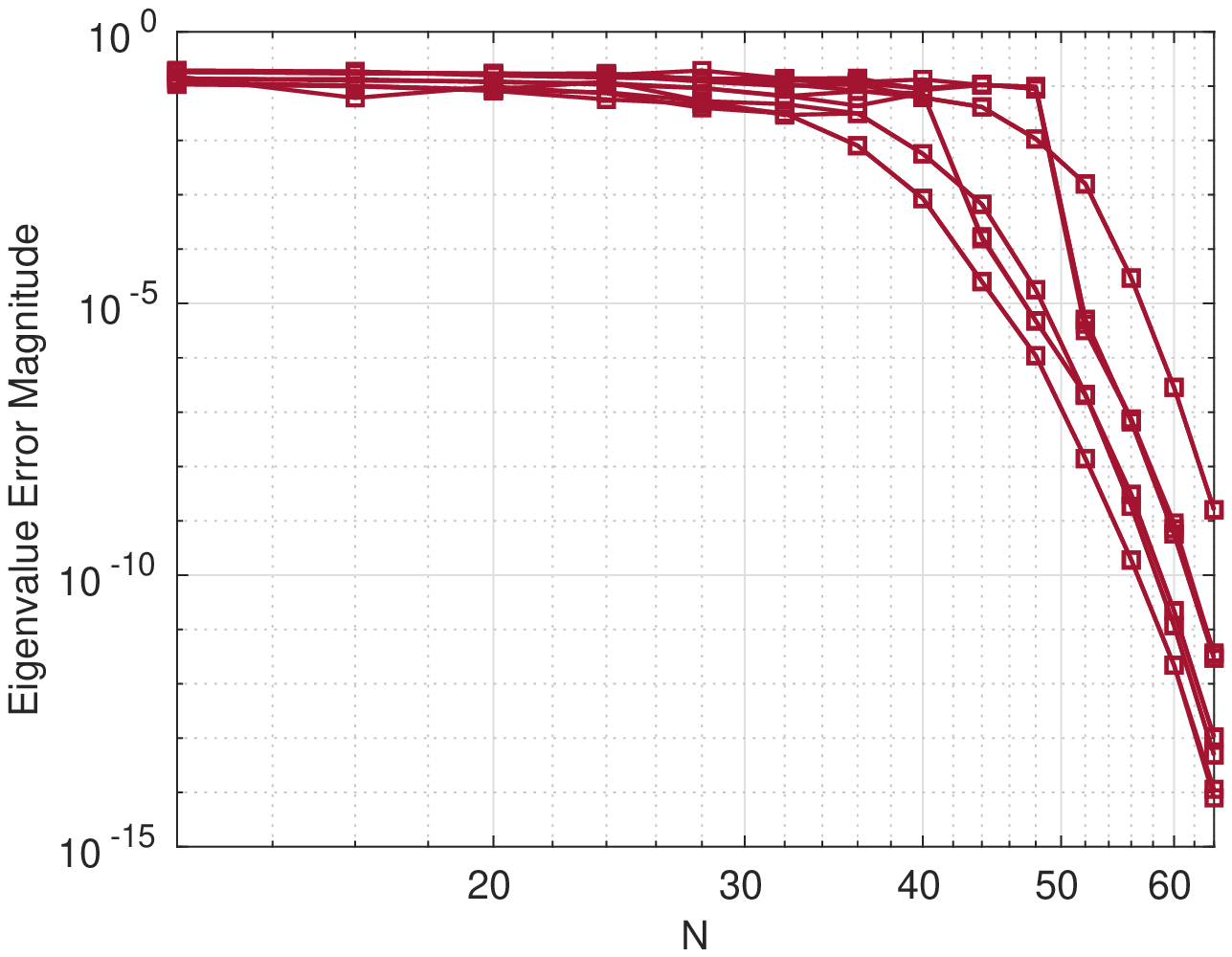}
\caption{Spectral convergence of the first 12 eigenvalues for $\nu=0.01$ (left) and $\nu=0.001$ (right).\label{fig:eigconvergence}}
\end{figure}

To show convergence of the algorithm, we track the error in the first 12 eigenvalues $\left\{ \lambda_{N,j}^{(\nu)} \right\}_{j=1}^{12}$ for different values of $N = 12, 16, 20, \dots, 64$ with respect to a more refined set obtained with $N = 80$. Here, we fix $\omega_0 = 0$, $r = 0.5$, $\beta(x) = \cos(x_1)+\sin(x_2)$ and take $\nu = 0.01$ and $\nu = 0.001$. We see in Figure \ref{fig:eigconvergence} that the error $\left| \lambda_{N,j}^{(\nu)} - \lambda_{80,j}^{(\nu)} \right|$ decays spectrally fast in $N$, in agreement with the smoothness of the eigenfunctions that is expected from the elliptic perturbation to the operator $P(x,D)$.

\section{Singular behaviour of the long-term evolution} \label{sec:singular_beh}

Even in the presence of infinitely smooth data, the development of attractors will produce singularities in $u$. Using microlocal analysis, results \cite{cs20,dz19} have been able to quantify this regularity. Moreover, they reveal a rather distributional character of the solution. The following result is a consequence of \cite[Theorem 3.1]{cs20} and \cite[Theorem]{dz19}.

\begin{theorem} \label{th:uinf}
Let $f \in C^{\infty}(\T^2)$, $\omega_0$ not an eigenvalue of $P$. Then, the solution to \eqref{eq:evol_problem} can be uniquely decomposed as
\begin{equation}
u(t) = e^{-i\omega_0 t} u_\infty + b(t) + \epsilon(t),
\end{equation}
where
\begin{enumerate}
\item $u_\infty := \lim_{\varepsilon \to 0}(P - \omega_0 -i\varepsilon)^{-1}f$ belongs to $\Hs{s}$ for any $s<-1/2$ and is not in $\Lt$ except if it vanishes,
\item $b$ is a bounded function with values in $\Lt$,
\item $\epsilon$ vanishes as $t \to \infty$ in the $\Hs{s}$-norm, for any $s<-1/2$.
\end{enumerate}
Moreover, the energy $\norm{u(t)}_0^2$ grows linearly except if $u_\infty$ vanishes.
\end{theorem}

While this result is stated in \cite[Theorem 3.1]{cs20} for a general bounded, self-adjoint, pseudo-differential operator of degree 0 that satisfies certain dynamical assumptions, the structure of the operator defined in \eqref{eq:def_P} has been considered in \cite{dz19} as a feasible choice for this problem. Furthermore, \cite{dz19} also proves a similar version of Theorem \ref{th:uinf}, but the proof relies instead in standard radial estimates (cf. \cite[\textsection E.4]{dzbook}). 

\subsection{Regularity of the computed approximations}

To illustrate how some of the statements in Theorem \ref{th:uinf} manisfest in the numerical experiments (and in particular, in the regularity of an approximation $u_N$), we have a look at how their corresponding Fourier coefficients decay. In two dimensions, the analysis can be done in a radial fashion (see, e.g., \cite{cgss2015,rnt12}).

\begin{definition}[Radial Energy Density (RED)] \label{def:red}
Let $N \geq 4$ be an even integer and $u_N$ be a discrete function defined on the grid $\Tau_N$. For $s \in \R$, we define the radial energy density (RED) $E_s$ of a discrete function $u_N$ as:
\begin{equation} \label{eq:def_red}
E_s[u_N](R) := \dfrac{1}{N^2} \sum_{k \in A_{R} \cap \Z^2} \ip{k}^{2s} \, |\widehat{u}_N(k)|^2, \quad R = 2,4,6, \dots, \frac{N}{2},
\end{equation}
where $A_{R} := \Big\{ x\in \R^2 : R-2 \leq |x| < R \Big\}$ is the $R$-th annulus of width $2$ in $\R^2$.
\end{definition}

\begin{remark} \label{re:red}
From the previous definition, we readily see that
\[
h^2 \sum_{l=1}^{N/4} E_s[u_N](2l) \leq \dfrac{h^2}{N^2} \sum_{k_1=-N/2}^{N/2-1}\sum_{k_2=-N/2}^{N/2-1} \ip{k}^{2s} \, |\widehat{u}_N(k)|^2  = \nnorm{u_N}_s^2,
\]
and therefore, for large $N$, if $u_N \in \Hs{s}$ then the series $\sum_{l=1}^\infty E_s[u_N](2l)$ must converge.
\end{remark}

\begin{figure}[t]
\centering
\includegraphics[width=0.5\textwidth]{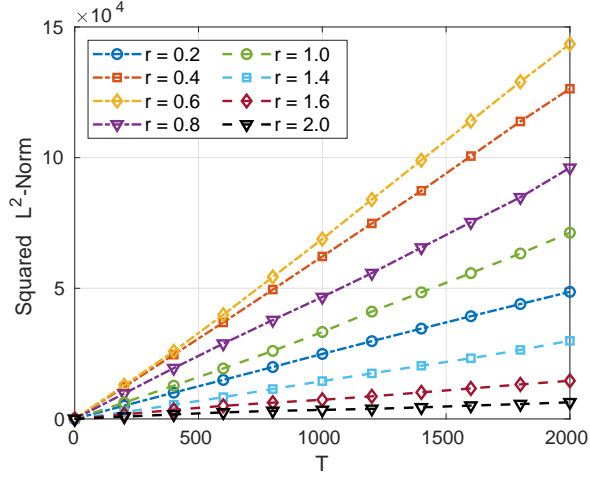}
\caption{Evolution of the squared $\rm L^2$-norm for several values of $r$ and $\beta(x) = \cos(x_1)+\sin(x_2)$.\label{fig:sql2norm}}
\end{figure}

\begin{figure}[t]
\centering
\includegraphics[width=\textwidth]{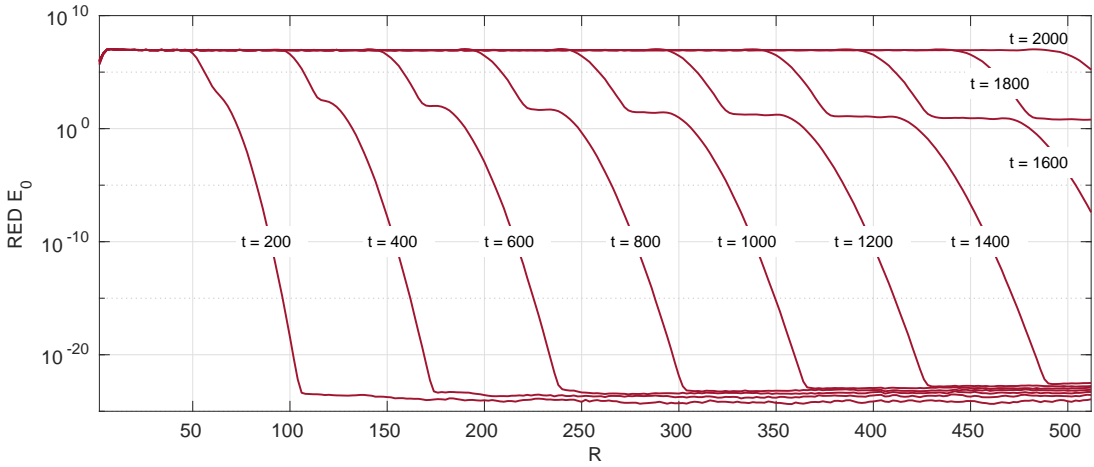}
\caption{Radial energy density $E_0$ at several times ($r = 0.25$, $\beta(x) = \cos(x_1)+\sin(x_2)$).\label{fig:sol_red0}}
\end{figure}

\begin{figure}
\centering
\includegraphics[width=\textwidth]{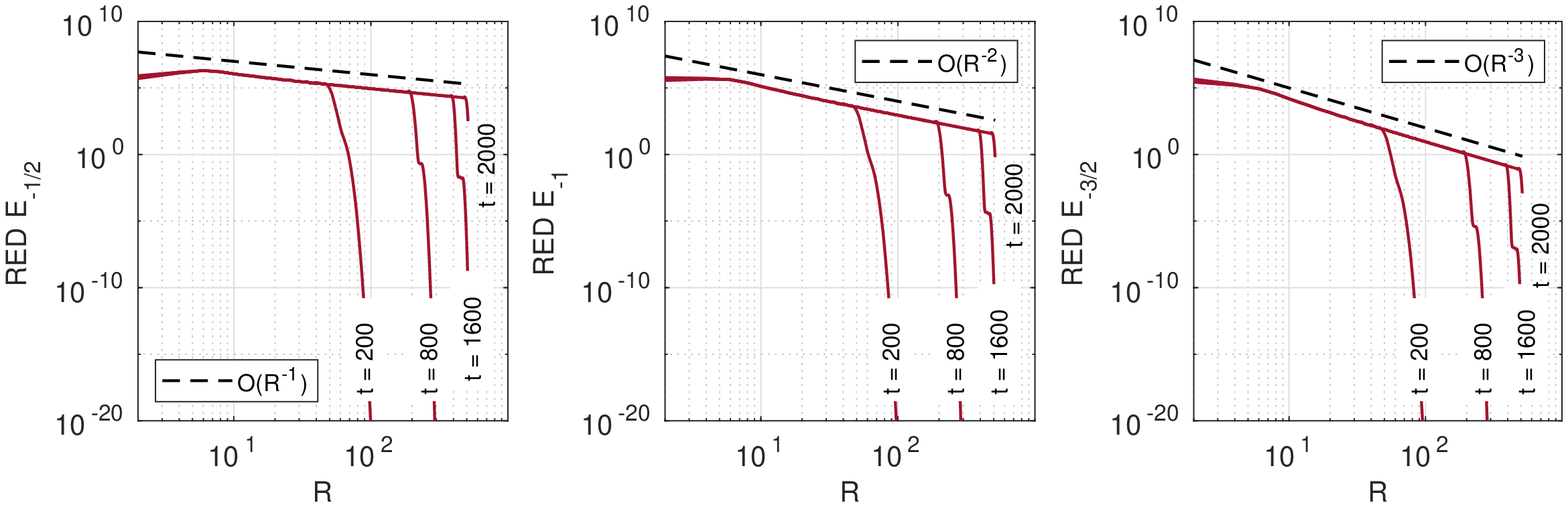}
\caption{Radial energy density for different values of $s$. This suggests that $E_s(u_N)(R)$ decays faster than $R^{2s}$ for any $s \leq 0$.\label{fig:sol_reds}}
\end{figure}

First, we show in Figure \ref{fig:sql2norm} the linearity in the evolution of $\nnorm{u_N(t)}_0^2$ for different values of $r$ (the source term is again a centered Gaussian similar to the one used in Figure \ref{fig:attractors}). 
It becomes more evident that $u_\infty \notin \Lt$ when we look at the RED $E_0$ (cf. \eqref{eq:def_red}) at different times, as shown in Figure \ref{fig:sol_red0}: while the RED quickly drops to below machine epsilon at the beginning of the simulation, it is not the case as the end time $T$ increases. Finally, we can have a look at how fast the RED decays for several choices of $s$. 
First, for $s=-1/2$, we observe in Figure \ref{fig:sol_reds} (left) that the RED decays as $R^{-1}$, which is slow enough to say that $u_\infty \notin \Hs{-1/2}$ (per Remark \ref{re:red}, since adding all points in the $E_{-1/2}[u_N]$ curve would resemble the harmonic series). However, as soon as we take $s < -1/2$, the RED appears to decay as $R^{2s}$, as shown as in Figure \ref{fig:sol_reds} (center and right), which suggests that $u_\infty \in \Hs{s}$ for $s < -1/2$.

\subsection{Attractors and energy manifolds} \label{sec:attractors_longterm}

As we mentioned in Section \ref{sec:background}, understanding the flow of the (rescaled) Hamiltonian vector field over the energy manifold is fundamental to understand the dynamics of the system. The long-term evolution described by Theorem \ref{th:uinf} was first given in the work by Colin de Verdi\`ere \& Saint-Raymond \cite{cs20}. There, two of the main assumptions for this flow are:
\begin{enumerate}
\item The flow of $X$ on $\Sigma$ is Morse-Smale with no fixed points. By the Poincar\'e-Bendixon theorem (cf., e.g., \cite{nz99}), this forces $\Sigma$ to be a finite union of tori.
\item The energy surface $\Sigma$ covers $\T^2$.
\end{enumerate}  
Although the latter has been relaxed in later works (see \cite{cdv18, dz19}), we believe that a loss of this property may change the behaviour of the solution in a noticeable way. After all, the attractors form a set $\Lambda_0^+$ whose projection by $\kappa$ lives on the energy manifold $\Sigma \subset \partial \tbs \T^2$.

Note that $\Sigma$ is an orientable surface given by the zeros of the principal symbol of $P(x,D)$ in the orientable 3-manifold $\partial \tbs \T^2$. This allows us to plot these manifolds for the operator in study. Similar to Section \ref{sec:example_cos}, for $\omega_0 = 0$, the principal symbol of $P(x,D)$ is given by
\[
\bar{p}(x,\xi) = |\xi|^{-1} \xi_2 - r\beta(x), \quad (x,\xi) \in \tbs{\T^2} \backslash \{0\}
\]
Then, for $(x,\xi) \in \Sigma \subset \partial \tbs \T^2$, we can parametrize $\xi$ as $(s \cos(\eta),s \sin(\eta)$ for $\eta \in [-\pi,\pi)$ and some $s > 0$. Thus, the energy manifold $\Sigma$ can be characterized as\begin{equation}
\Sigma  = \Big\{ (x_1,x_2,\eta) \in \T^3 : \ r\, \beta(x_1,x_2) = \sin(\eta)  \Big\},
\end{equation} 
where $\T^3$ is the standard 3-torus. In particular, relating to the second assumption mentioned at the beginning of this subsection, $\Sigma$ will not cover $\T^2$ if and only if there exists $x=(x_1,x_2) \in \T^2$ such that for any $\eta \in \mathbb{S}^1$, the equation $r\beta(x) = \sin(\eta)$ does not have a solution.

Let us consider some choices of data for which we can compare the long-term evolution of the solution (where the attractors are fully developed) and the mentioned energy manifolds:

\begin{description}
\item[Test 1:] $r=0.5$, $\beta(x) = \cos(x_1)$,
%\item[Test 2:] $r=0.5$, $\beta(x) = \cos(x_1)+\sin(x_2)$ and $\nu \in [1.3,9.3]\cdot 10^{-3}$.
\item[Test 2:] $r=0.45$, $\beta(x) = \cos(x_1-2x_2)+\sin(2x_2)$,
\item[Test 3:] $r=0.55$, $\beta(x) = \cos(x_1-2x_2)+\sin(2x_2)$.
\end{description}

In all cases we take $\omega_0 = 0$. The respective energy manifolds read:
\begin{align}
\label{eq:def_s1}\Sigma_1 &:= \Big\{ (x_1,x_2,\eta) \in \T^3 \ : \ 0.5 \cos(x_1) = \sin(\eta) \ \Big\}, \\
\label{eq:def_s2}\Sigma_2 &:= \Big\{ (x_1,x_2,\eta) \in \T^3  \ : \ 0.45\Big(\cos(x_1-2x_2)+\sin(2x_2) \Big) = \sin(\eta) \ \Big\}, \\
\label{eq:def_s3}\Sigma_3 &:= \Big\{ (x_1,x_2,\eta) \in \T^3  \ : \ 0.55\Big(\cos(x_1-2x_2)+\sin(2x_2) \Big) = \sin(\eta) \ \Big\}.
\end{align}

\begin{figure}
\centering
\includegraphics[width=0.99\textwidth]{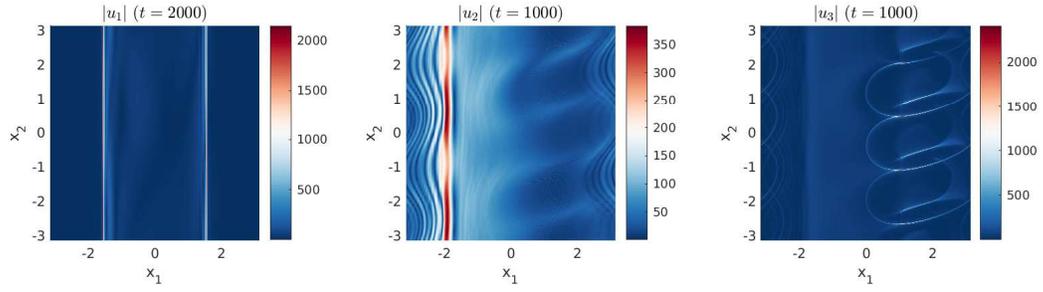}
\caption{Long-term evolution for Tests 1, 2 3 (left to right). Evolution computed using $\omega_0 = 0$ and a centred Gaussian as a source term.\label{fig:lte}}
\end{figure}

\begin{figure}
\includegraphics[width=0.99\textwidth]{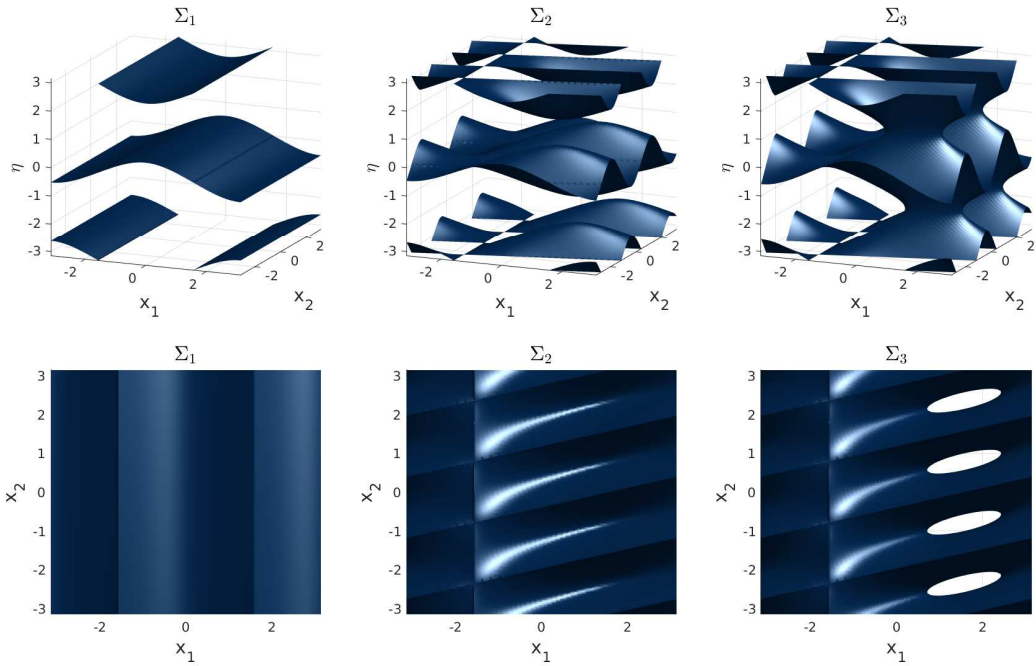}
\caption{Energy manifolds $\Sigma_j$, $j \in \{1,2,3\}$ viewed from the side (top row), and from the top (bottom row). This last view reveals holes in $\Sigma_3$, but not in $\Sigma_1$ and $\Sigma_2$. \label{fig:manifolds}}
\end{figure}

Notice that Test 1 corresponds to the data used to introduce this work (see Figure \ref{fig:attractors}). Also note that the difference in Test 2 and 3 is only a slight increase in the parameter $r$. For each one of the tests, the solution to the evolution \eqref{eq:evol_problem} containing the attractors is shown in Figure \ref{fig:lte} and the energy manifolds $\Sigma_j$ are shown in Figure \ref{fig:manifolds}. Here we see how these manifolds shape the attractors. In particular, the increase from $r=0.45$ in Test 2 to $r=0.55$ in Test 3 causes the energy manifold $\Sigma_3$ to not cover the 2-torus, but attractors still develop in Test 3.  This not only shows that the second assumption is not needed but also how a small variation in the parameters can completely change the shape of the attractors.

\section{Spectra of vanishing-viscosity operators} \label{sec:spectra}

% Discussion of the spectral properties of original operator, the viscous operator, and how you ORDERED the eigenvalues. 

For rotating fluids, it is known that the existence of internal wave attractors is related to the spectral properties of the underlying differential operator (cf. \cite{Ralston1973,RGV01}). This motivates the study of the eigenvalue problem related to \eqref{eq:evol_problem}:
\begin{equation} \label{eq:ev_problem_P}
(P(x,D)-\omega_0)u(x) = \lambda u(x), \quad x \in \T^2,
\end{equation}
with $\omega_0 \geq 0$ and $P$ as defined in \eqref{eq:def_P}. Notice that, without loss of generality, we can focus on the case $\omega_0 = 0$, since $\spec{}{P-\omega_0} = \spec{}{P}-\omega_0$ (by the spectral mapping theorem, cf. \cite[Theorem 7.1.9]{bogachev}), so the eigenvalues only get shifted by $\omega_0$. Also, since $P$ is a bounded self-adjoint operator, the spectrum $\spec{}{P}$ lies entirely on the real axis (cf., e.g., \cite[Corollary 7.2.5]{bogachev}). We can obtain a more precise characterization of the essential spectrum using \cite[Lemma 2.1]{cdv18}:
\begin{equation} \label{eq:ess}
\spec{ess}{P} \subset \left[ -1-r\max_{\T^2} \beta \, , \, 1-r\min_{\T^2} \beta \right].
\end{equation}
Here, the essential spectrum has to be understood in the sense of Pelinovsky \cite{pelinovsky2011}, that is, $\spec{ess}{P} := \spec{}{P} \backslash \spec{d}{P}$, where the discrete spectrum $\spec{d}{P}$ is the set of all eigenvalues of $P$ with finite (algebraic) multiplicity and which are isolated points of $\spec{}{P}$. 

\begin{remark}
According to \cite{cs20}, the presence of a continuous spectrum in $P(x,D)$ as in \eqref{eq:ess} and the assumptions on the flow of the Hamiltonian field as in Section \ref{sec:attractors_longterm} ensure the generation of attractors in the evolution of $P(x,D)$.
\end{remark}

It was shown in \cite{dz19} that there exists a finite number of eigenvalues. Moreover, their associated eigenfunctions that are analytic (cf. \cite{Wang2020}). The following result summarizes these findings.

\begin{lemma}[{\cite[Lemma 3.2]{dz19}}] \label{th:spec_P} Let $\omega_0 = 0$. There exists $\delta$ sufficiently small such that { the cardinality of $\spec{pp}{P} \cap [-\delta,\delta]$ is finite}. Furthermore, if $Pu=\lambda u$ for $u \in \Lt$ and $|\lambda|\leq \delta$, then $u \in C^\infty(\T^2)$. 
\end{lemma}

Here, the {\it pure point spectrum} $\spec{pp}{P}$ is the set of all eigenvalues, both embedded and isolated (a precise definition can be made using the spectral measure, cf. \cite[\S 4.3]{amrein}). This opens the possibility of having eigenvalues that are embedded in the essential spectrum (i.e., at a zero distance from the continuous spectrum $\spec{c}{P}$) given by \eqref{eq:ess}, especially when $[-\delta,\delta] \subset \spec{ess}{P}$. In this case, a straightforward discretization of \eqref{eq:ev_problem_P} would not differentiate between points in $\spec{pp}{P}$ or $\spec{ess}{P}$, so we need to take a different approach such that these elements become distinguishable.

\subsection{Elliptic perturbation} \label{sec:elliptic_perturbation}

Physics literature (such as \cite{RGV01}) suggests that, for rotating fluids, we can get more information about these eigenvalues by perturbing the equation with a small viscous term. Applying this idea to stratified fluids means that the zeroth-order operator $P$ converts into the second-order operator $P+i\nu\Delta$ (where $\Delta$ is the standard Laplacian), however, these two operators are completely different. The first one has a combination of continuous spectrum and embedded eigenvalues, while the second one has a purely discrete spectrum that is much easier to compute numerically. Recent result \cite{gz19} focuses on the study of the limit of $P+i\nu\Delta$ as $\nu \to 0^+$. In particular, the following result, which justifies the discretization of the eigenvalue problem \eqref{eq:ev_problem}, is due to \cite[Theorem 1]{gz19}.

\begin{theorem} \label{th:viscous}
Consider the operator $P+i\nu\Delta$, with $P$ as given in \eqref{eq:def_P}. Then, there exists an open neighbourhood $U$ of $0$ in $\C$, and a set
\[
\RR(P) \subset \{ z \in \C : \Im(z) \, \leq 0 \} \cap U,
\]
such that for every set $K$ compactly contained in $U$, $\RR(P) \cap K$ is discrete and
\[
\spec{pp}{P+i\nu\Delta} \cap U \longrightarrow \RR(P) \quad \text{as }\nu \to 0^+,
\]
uniformly on $K$. Furthermore,
\[
\RR(P) \cap \R = \spec{pp}{P} \cap U.
\]
\end{theorem}

Here, the set $\RR(P)$ is known as the set of \textit{resonances} of $P$. When  restricted to $U$, this set is made of all eigenvalues of the limiting operator $\lim_{\nu \to 0^+} P+i\nu\Delta$ that are contained in $U$. Because this operator is not self-adjoint, not all eigenvalues will lay on the real axis, but those who do, will be precisely the embedded eigenvalues of $P$.  

The uniformity in the convergence stated in Theorem \ref{th:viscous} takes a key role in numerical approximations. It tells us that every eigenvalue of $P$ that lies in the neighbourhood $U$ can be approximated by viscous eigenvalues. More precisely, let us write $\RR(P) = \left\{ \lambda_j \right\}_{j=1}^N$, where $N=\infty$ is allowed (recall that $\RR(P)$ not only contains some of the embedded eigenvalues of $P$ but also elements from the resolvent set $\rho(P) := \C \backslash \spec{}{P}$). The previous theorem suggests that for $\spec{pp}{P+i\nu\Delta} = \left\{\lambda_j^{(\nu)} \right\}_{j=1}^\infty$, we have (after suitable reordering)
\begin{equation} \label{eq:conv_evals}
\lambda_j^{(\nu)} \longrightarrow \lambda_j \quad \text{as }\nu \to 0^+,
\end{equation}
uniformly on compact sets and with agreement of multiplicities. This means that, by tracking the limit $\nu \to 0^+$, we should be able to find some of the embedded eigenvalues (and eigenfunctions) of $P$.

One factor that must be considered though is that, while the embedded eigenfunctions $u$ of $P$ are analytic, their approximations given by the eigenmodes $u^{(\nu)}$ of $P+i\nu\Delta$ might have poor regularity (for fixed $\nu$). Indeed, \cite[Theorem 2]{gz19} affirms that there exists $\widetilde{\delta} > 0$ such that the Hilbert space $\mathcal{X}$ in which the eigenmodes $u^{(\nu)}$ live satisfies:
\begin{equation} \label{eq:hyperfunctions}
\mathscr{A}_{\widetilde{\delta}} \subset \mathcal{X} \subset \mathscr{A}_{-\widetilde{\delta}},
\end{equation}
where for any $s \geq 0 $, $\mathscr{A}_{-s}$ is the dual of the space $\mathscr{A}_s$, defined as
\[
\mathscr{A}_s:= \left\{ u \in \Lt : \sum_{\xi \in \Z^2} |\widehat{u}(\xi)|^2 e^{4|\xi|s} < \infty \right\}.
\]

\subsection{Ordering of eigenvalues} \label{sec:ordering}

As $\nu \to 0^+$, the eigenvalues $\lnu$ draw defined curves in the complex plane that are smooth \cite{Wang2020}. To track these trajectories, the computed eigenvalues must be sorted in an appropriate way. Thus, after computing the closest eigenvalues to $\omega_0$ (as explained in Section \ref{sec:3_ev_problem}), we choose to order them in two steps:
\begin{enumerate}
\item Sort using a ``magnitude-then-phase'' approach, that is, first the eigenvalues are ordered in increasing magnitude, and if two eigenvalues have the same magnitude, the one with smallest phase (in the interval $(-\pi,\pi]$ goes first. Then,
\item Move all eigenvalues with nonnegative real part to the top of the column vector containing the requested eigenvalues.
\end{enumerate} 

Unfortunately, the shape of these curves is highly problem-dependent, and this sorting procedure may not be useful if, for instance, there is an eigenvalue with zero real part (since numerically this 0 could manifest as, say, $\pm 10^{-16}$). In this case, the sorting must be undone, and proceed with a different approach.

\subsection{Resonances near the origin} \label{sec:resonances_near_0}

\begin{figure}[t]
\centering
\includegraphics[width=\textwidth]{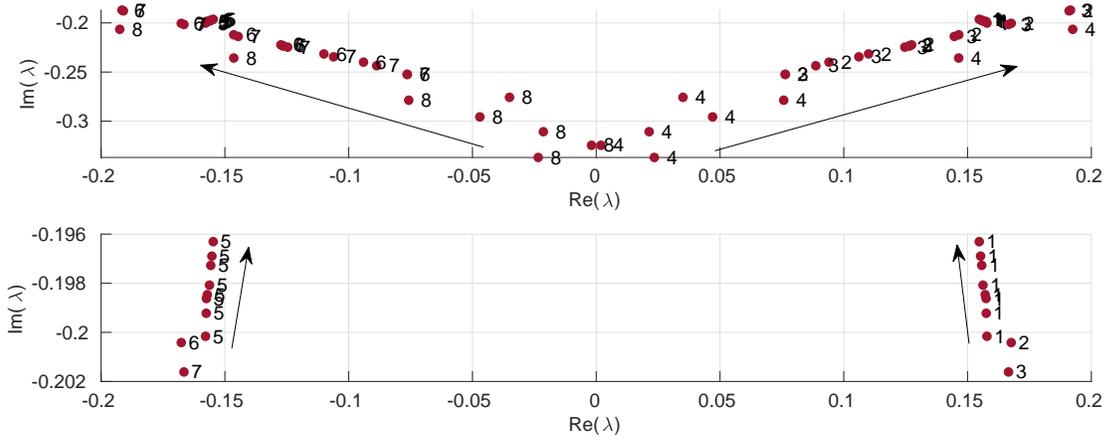}
\caption{Evolution of the first 8 eigenvalues in Test 1 ($r=0.5$, $\beta(x) = \cos(x_1)$) as $\nu$ decreases from $9.3 \cdot 10^{-3}$ to $2.3 \cdot 10^{-3}$. While some of them move in an oblique direction away from 0 (top), a closer look near the points $\pm 0.15 - 0.2i$ (bottom) shows that the first and fifth eigenvalues are moving slowly upwards.\label{fig:t1_evals}}
\end{figure}

\begin{figure}
\centering
	\includegraphics[width=\textwidth]{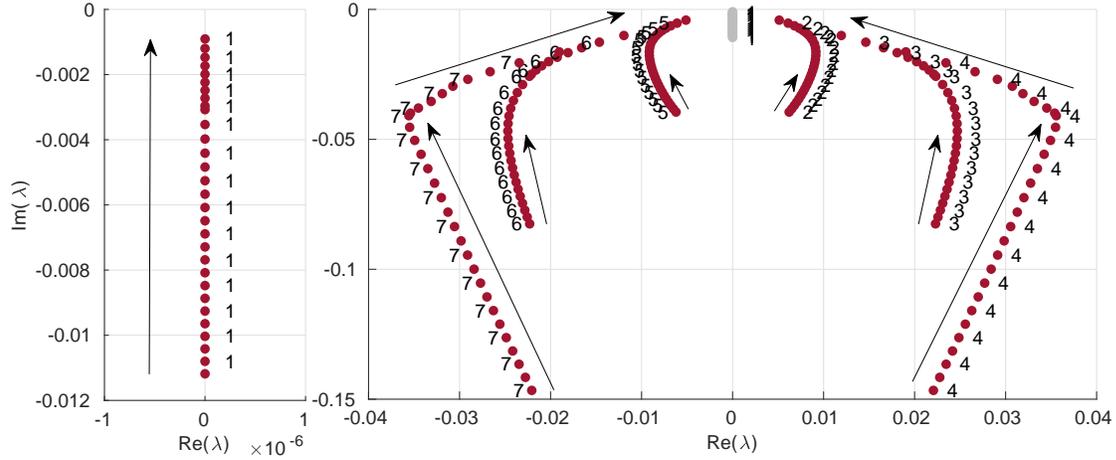}
\caption{Evolution of the first eigenvalue (left) and second to seventh eigenvalues (right) for Test 2 ($r=0.45$, $\beta(x) = \cos(x_1-2x_2)+\sin(2x_2)$) when the viscosity decreases from $10^{-2}$ to $3 \cdot 10^{-4}$.\label{fig:t2_evals}}
\end{figure}

\begin{figure}
\centering
	\includegraphics[width=\textwidth]{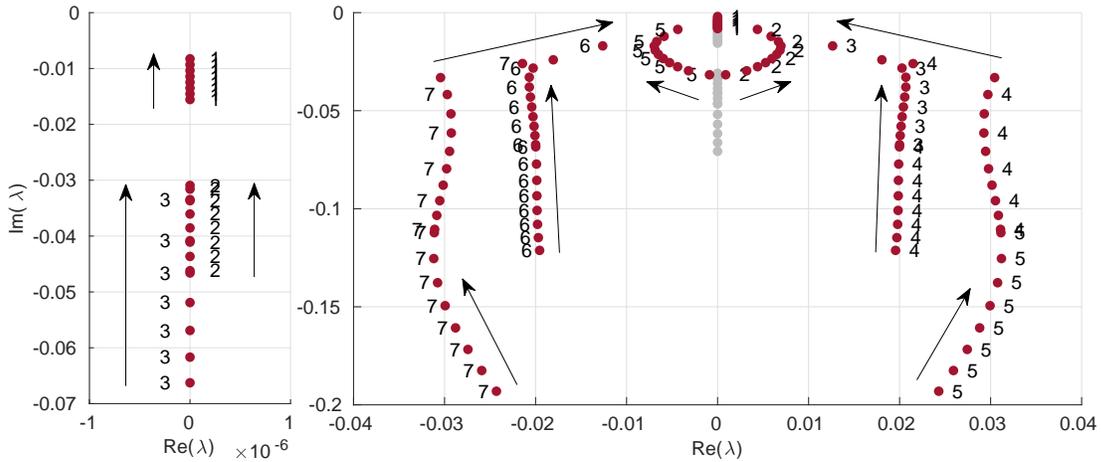}
\caption{Evolution of the first seven eigenvalues for Test 3 ($r=0.55$, $\beta(x) = \cos(x_1-2x_2)+\sin(2x_2)$). Left: first three eigenvalues when the viscosity decreases from $10^{-2}$ to $4.4 \cdot 10^{-3}$ (below this viscosity, these eigenvalues become the first, second and fifth eigenvalues, respectively).  Right: all seven eigenvalues as $\nu$ decreases from $10^{-2}$ to $3 \cdot 10^{-4}$.\label{fig:t3_evals}}
\end{figure}

Let us compute some of these eigenvalues and track their trajectories as $\nu \to 0^+$. We consider the same list of Tests as in Section \ref{sec:attractors_longterm}, but with various small viscosities:
\begin{description}
\item[Test 1:] $r=0.5$, $\beta(x) = \cos(x_1)$, and $\nu \in [2.3,9.3]\cdot 10^{-3}$,
\item[Test 2:] $r=0.45$, $\beta(x) = \cos(x_1-2x_2)+\sin(2x_2)$, and $\nu \in [3\cdot 10^{-4}, 10^{-2}]$,
\item[Test 3:] $r=0.55$, $\beta(x) = \cos(x_1-2x_2)+\sin(2x_2)$, and $\nu \in [3\cdot 10^{-4}, 10^{-2}]$.
\end{description}
For Test 1, we compute the first 8 eigenvalues of the operator $P+i\nu\Delta$, whereas for Test 2 and 3 we compute the first 7 eigenvalues. We portray these results in Figures \ref{fig:t1_evals}, \ref{fig:t2_evals}, and \ref{fig:t3_evals}. In all the experiments, the eigenvalue problem \eqref{eq:ev_problem} is discretized using a mesh $\Tau_N$ with $N=64$ per direction.

Overall, the eigenvalues are located in the lower half of the complex plane, and these move upwards toward the real axis as $\nu \to 0^+$. This is expected since $i\nu\Delta$ is a second-order differential operator with a purely complex spectrum that lies on the negative part of the imaginary axis. However, as pointed out in Section \ref{sec:elliptic_perturbation}, in the limit there might be some eigenvalues that will stay below the real axis (recall that $P+i\nu\Delta$ is not a self-adjoint operator and that $\spec{}{P+i\nu\Delta} \centernot \longrightarrow \spec{}{P}$ in its entirety as $\nu \to 0^+$). 

Notice that the ordering of these eigenvalues works well for Test 1 in the sense that we are able to track important trajectories (such as ones drawn by the first and fifth eigenvalues). Furthermore, we confirm that the trajectories are smooth. We also see that there is an observed symmetry with respect to the imaginary axis.  Figure \ref{fig:t1_evals} suggest that if $\lambda^{(\nu)}$ is an eigenvalue, then so is $-\overline{\lambda^{(\nu)}}$. However, the situation is different in Tests 2 and 3. 

First, we see in Figure \ref{fig:t2_evals} that $\lambda_1^{(\nu)}$ is an simple eigenvalue moving along the real axis. Moreover, in Test 3 (see Figure \ref{fig:t3_evals}), $\lambda_2^{(\nu)}$ and $\lambda_5^{(\nu)}$ appear in pairs (as before) but only below $\nu = 4.4 \cdot 10^{-3}$. Therefore, for these particular cases, we have ordered the eigenvalues in a ``magnitude-then-phase'' way first (described in Section \ref{sec:ordering}), and then used the following order: first those eigenvalues $\lambda$ with $|\Re(\lambda)| <= 0.5\cdot10^{-3}$, then those with $\Re(\lambda)>0.5\cdot10^{-3}$ and finally those with $\Re(\lambda) < -0.5\cdot10^{-3}$.  

One additional thing to notice in Tests 2 and 3 is how all seven computed eigenvalues move toward a neighbourhood of 0 as $\nu \to 0^+$, in comparison to Test 1 where some eigenvalues appear to be moving toward the real axis, but not toward 0. We conjecture that the difference in these cases is due to the presence of an eigenvalue at 0 with different multiplicities.

\subsection{Regularity of eigenmodes}

\begin{figure}
\centering
\begin{tabular}{|c|} \hline
Test 1: $\ip{D}^{-1} D_{x_2} + i\nu\Delta - 0.5\cos(x_1)$, $\nu \in [2.3 \cdot 10^{-3}, 9.3 \cdot 10^{-3}]$ \\\hline
\includegraphics[width=0.96\textwidth]{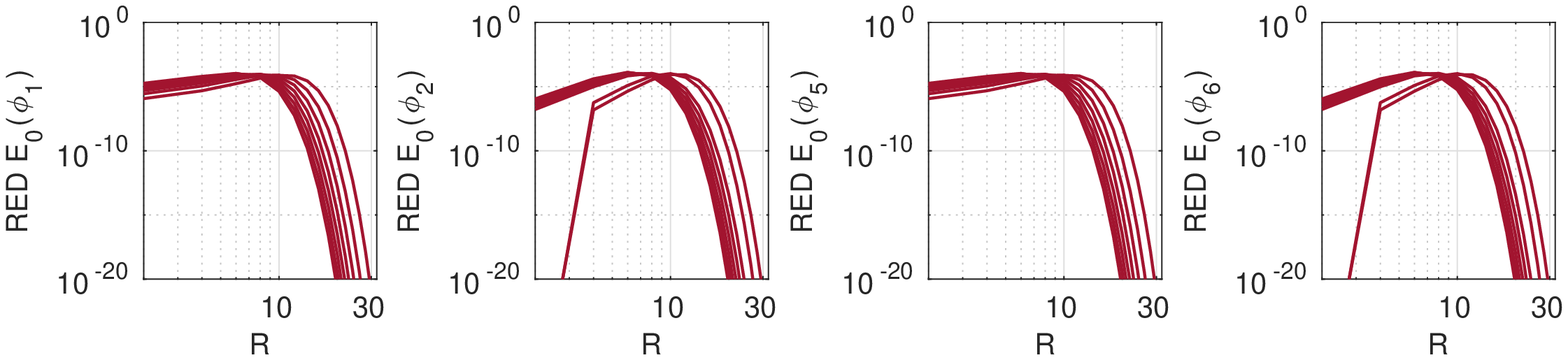} \\\hline
Test 2: $\ip{D}^{-1} D_{x_2} + i\nu\Delta - 0.45\left(\cos(x_1-2x_2) + \sin(2x_2)\right)$, $\nu \in [3 \cdot 10^{-4}, 2.4 \cdot 10^{-3}]$ \\\hline
\includegraphics[width=0.96\textwidth]{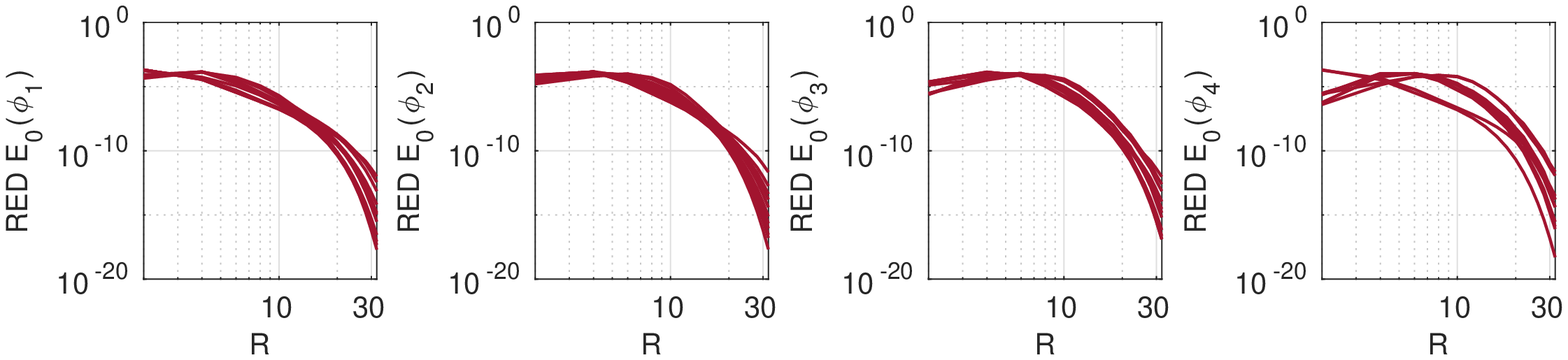} \\\hline
Test 3: $\ip{D}^{-1} D_{x_2} + i\nu\Delta - 0.55\left(\cos(x_1-2x_2) + \sin(2x_2)\right)$, $\nu \in [3 \cdot 10^{-4}, 4.3 \cdot 10^{-3}]$ \\\hline
\includegraphics[width=0.96\textwidth]{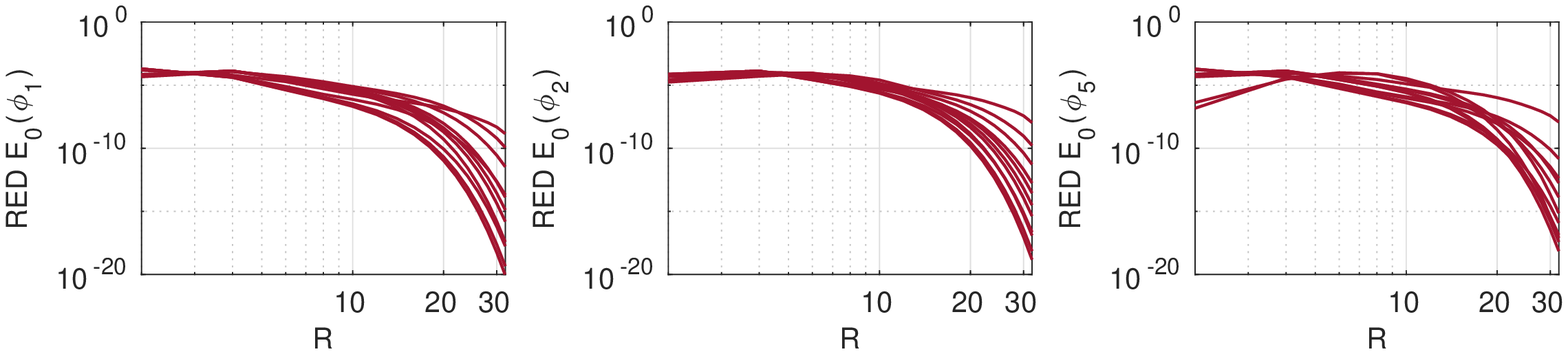} \\\hline
\end{tabular}
\caption{Radial energy density $E_0$ (log-log scale) of different viscous eigenfunctions $\phi_j$. Each curve represents a fixed value of $\nu$ within the chosen range. In general, as the viscosity decreases, the curves move to the right, which shows that the viscous eigenfunctions become less regular. \label{fig:red}}
\end{figure}

To get more information on the smoothness of some of the eigenmodes corresponding to $P+i\nu\Delta$, we have computed their radial energy density $E_0$ for the different viscosities considered in Figures  \ref{fig:t1_evals}, \ref{fig:t2_evals}, and \ref{fig:t3_evals}, with the focus mainly on those related to eigenvalues that appear to be moving to 0 (bearing in mind Theorem \ref{th:viscous}). These results are portrayed in Figure \ref{fig:red} for Tests 1, 2 and 3.

Overall, we observe in this figure that the Fourier coefficients decay faster when $\nu$ is large and slower when $\nu$ is small. The latter can be related to the poor regularity that the viscous approximations can have (see \eqref{eq:hyperfunctions}). We also notice how an increase in $r$ from 0.45 in Test 2 to 0.55 in Test 3 makes the Fourier coefficients of eigenfunctions decay slightly slower. This decrease in regularity may potentially be attributed to the fact that, for the choices in Test 3, the energy manifold $\Sigma_3$ does not cover $\T^2$ (see Section \ref{sec:attractors_longterm}).

\section{Low-viscosity eigenmodes and long-term behaviour} 

%HERE I WOULD LIKE TO WRITE ABOUT THE RELATIONSHIP BETWEEN EIGENMODES AND LTE, AS WELL AS TO DESCRIBE THE REGULARITY OF THESE EIGENMODES, SO THAT I CAN CONTRAST THEM WITH THE SINGULARITY OF THE LTE, ALL IN ONE SECTION.

We now present numerical evidence of a relationship between the internal wave attractors described in Section \ref{sec:singular_beh} and the spectra of zeroth-order operators discussed in Section \ref{sec:spectra}. More precisely, we will explore how the eigenfunctions of $P(x,D)$ can partially (and cannot completely) characterize the solution to the evolution problem \eqref{eq:evol_problem}.

While this characterization is complete for many elliptic operators, the fact that $P(x,D)$ has at most a finite number of eigenvalues (per Theorem \ref{th:spec_P}), their associated eigenfunctions cannot form a basis of $\Lt$. Consequently, it is neither immediate nor obvious that the eigenfunctions of this self-adjoint zeroth-order pseudo-differential operator could possibly describe the solution to the corresponding evolution problem. 

Let us have a look at some the viscous eigenmodes in \eqref{eq:ev_problem} corresponding to the eigenvalues closest to 0, with the viscosity taken as the smallest ones considered in Section \ref{sec:resonances_near_0}. We can then make the comparison with the long-term evolution of the solution to \eqref{eq:evol_problem}.

\begin{figure}
\centering
\begin{tabular}{|c|}\hline
Test 1: $\ip{D}^{-1} D_{x_2} + i\nu\Delta - 0.5\cos(x_1)$, $\nu = 2.3 \cdot 10^{-3}$ \\\hline
\includegraphics[width=0.96\textwidth]{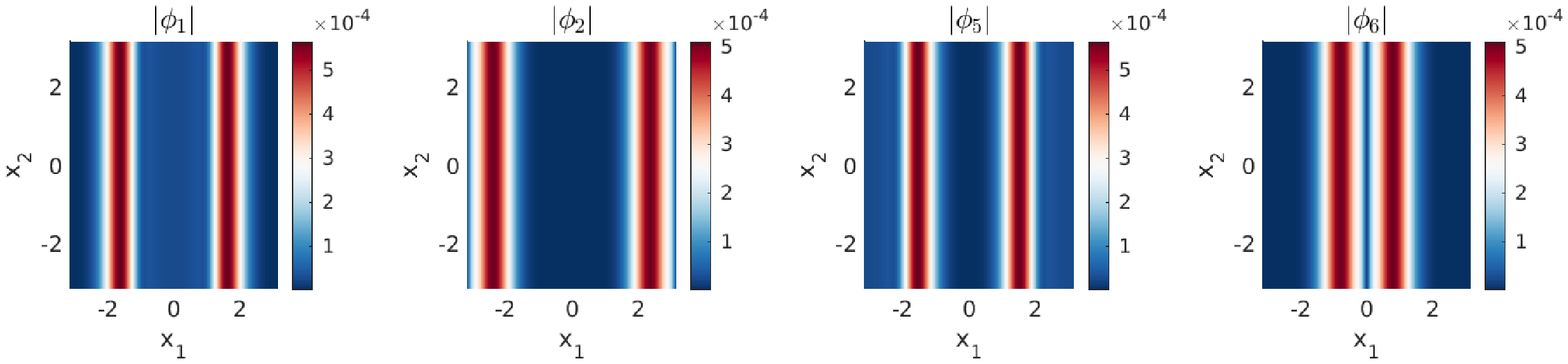} \\\hline
Test 2: $\ip{D}^{-1} D_{x_2} + i\nu\Delta - 0.45\left(\cos(x_1-2x_2) + \sin(2x_2)\right)$, $\nu = 3 \cdot 10^{-4}$ \\\hline
\includegraphics[width=0.96\textwidth]{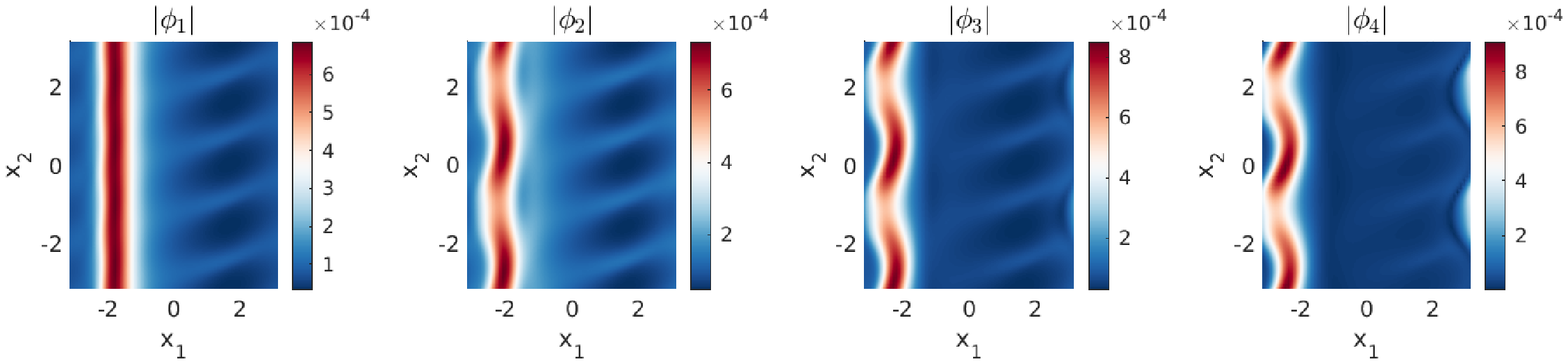} \\\hline
Test 3: $\ip{D}^{-1} D_{x_2} + i\nu\Delta - 0.55\left(\cos(x_1-2x_2) + \sin(2x_2)\right)$, $\nu = 3 \cdot 10^{-4}$ \\\hline
\includegraphics[width=0.96\textwidth]{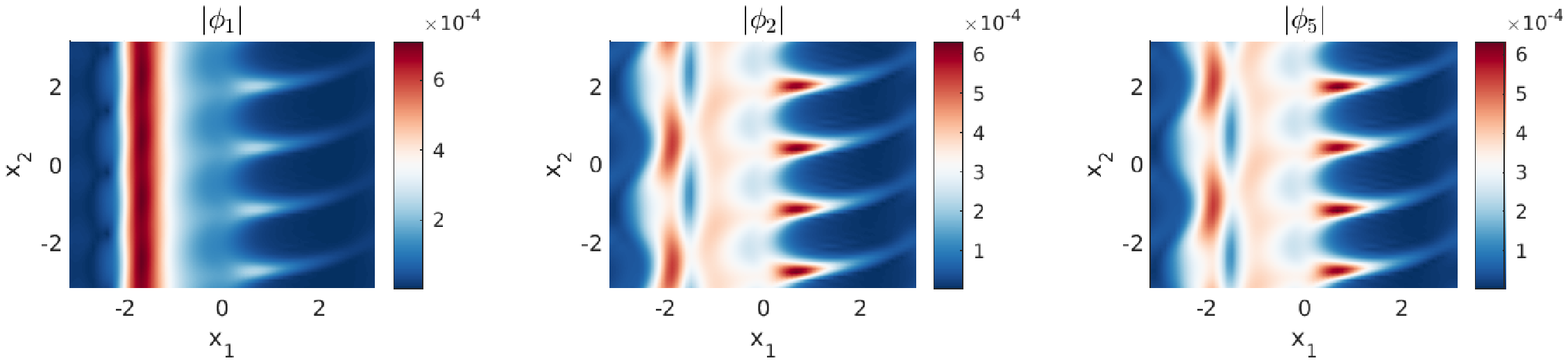} \\\hline
\end{tabular}
\caption{Magnitude of eigenmodes in real space $\phi_j$ for Tests 1, 2 and 3. The shape of some of the eigenmodes resembles that of the attractors in Figure \ref{fig:lte}. \label{fig:eigenmodes}}
\end{figure}

First, using the parameters from Test 1, we see in the first row of Figure \ref{fig:eigenmodes} how the first and fifth eigenmodes match the shape and location of the attractors in the transient solution (see Figure \ref{fig:attractors} or Figure \ref{fig:lte}). We can see similar situations using the parameters from Test 2 and 3 where different eigenmodes capture different parts of the attractors (compare Figure \ref{fig:eigenmodes} with Figure \ref{fig:lte}).

\begin{figure}
  \centering
  \includegraphics[width=0.99\textwidth]{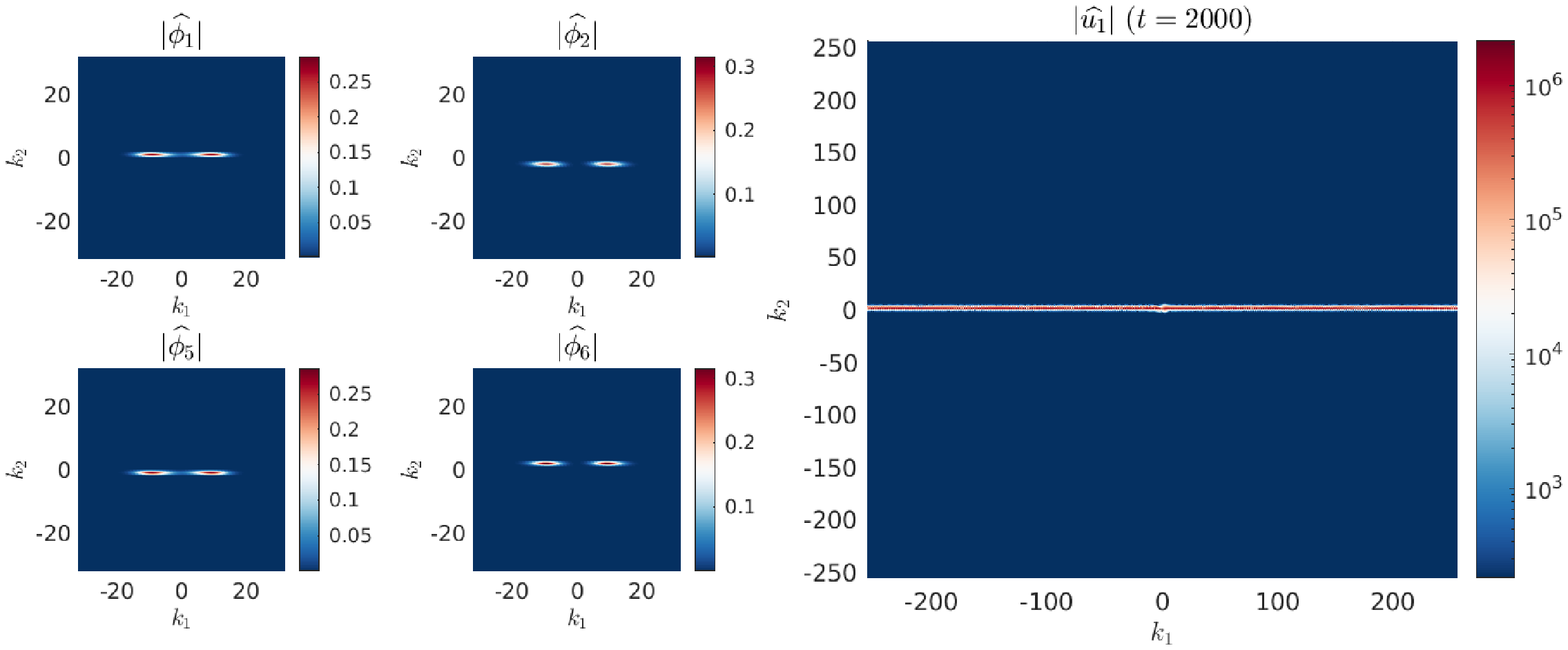}
  \caption{Left half: Magnitude of some eigenmodes in frequency space $\widehat{\phi_j}$ of $\ip{D}^{-1} D_{x_2} + i\nu\Delta - 0.5\cos(x_1)$ ($\nu = 2.3 \cdot 10^{-3}$). Right half: Long-term evolution in frequency space (this is the Fourier transform of Figure \ref{fig:lte}-left). \label{fig:eig_lte_1}}
  \end{figure}

\begin{figure}
  \centering
  \includegraphics[width=0.99\textwidth]{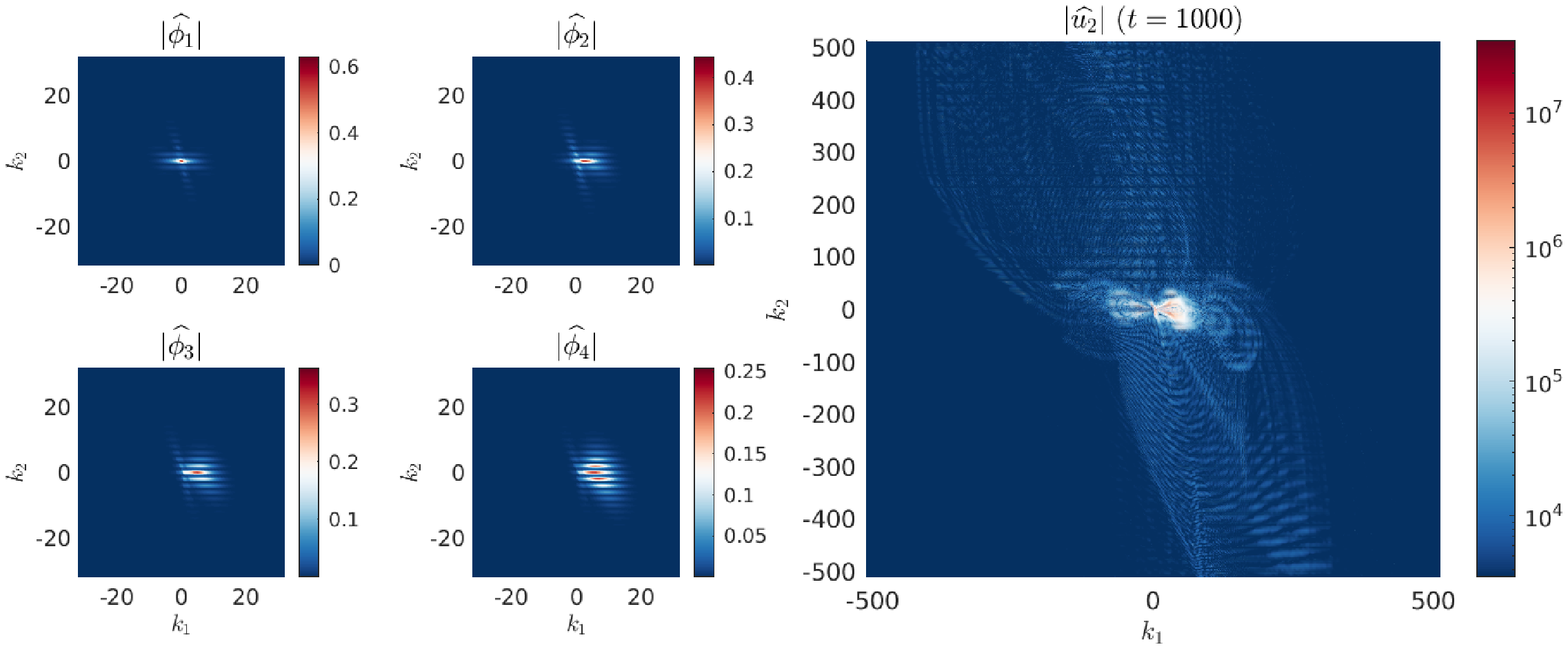}
  \caption{Left half: Magnitude of some eigenmodes in frequency space $\widehat{\phi_j}$ of $\ip{D}^{-1} D_{x_2} + i\nu\Delta - 0.45\left(\cos(x_1-2x_2) + \sin(2x_2)\right)$ ($\nu = 3 \cdot 10^{-4}$). Right half: Long-term evolution in frequency space (this is the Fourier transform of Figure \ref{fig:lte}-center). \label{fig:eig_lte_2}}
  \end{figure}

\begin{figure}
  \centering
  \includegraphics[width=0.99\textwidth]{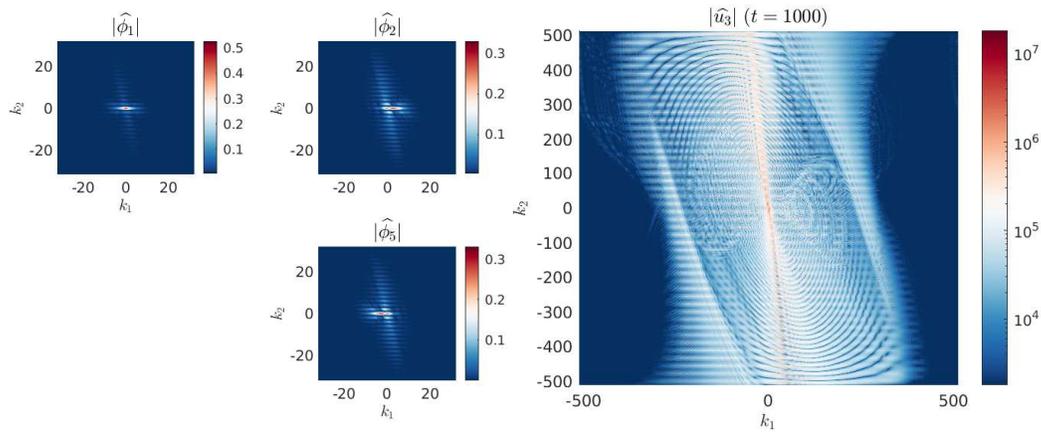}
  \caption{Left half: Magnitude of some eigenmodes in frequency space $\widehat{\phi_j}$ of $\ip{D}^{-1} D_{x_2} + i\nu\Delta - 0.55\left(\cos(x_1-2x_2) + \sin(2x_2)\right)$ ($\nu = 3 \cdot 10^{-4}$). Right half: Long-term evolution in frequency space (this is the Fourier transform of Figure \ref{fig:lte}-right). \label{fig:eig_lte_3}}
  \end{figure}

In turn, the magnitude of the Fourier coefficients of these modes seems also to provide some information about where the wave energy is concentrated in the transient solution. We portray this in Figures \ref{fig:eig_lte_1}, \ref{fig:eig_lte_2}, and \ref{fig:eig_lte_3}, respectively for Tests 1, 2, and 3. Additionally, by looking at these eigenmodes in Fourier space, we can corroborate the smoothness that the radial energy density $E_0$ suggests (discussed at the end of Section \ref{sec:spectra}). In particular, the eigenmodes in Test 1 and 2 seem to be compactly supported, whereas for Test 3, while the frequencies are more spread out, the highest amplitudes are still concentrated in the center of the spectrum.

Finally, we can also observe in Figures \ref{fig:eig_lte_1}-\ref{fig:eig_lte_3} how the numerics reveal the contrasting character between the embedded eigenmodes (represented by their viscous approximations in Figure \ref{fig:eigenmodes}, which we also portray in frequency space in \ref{fig:eig_lte_1}-\ref{fig:eig_lte_3}) and the long-term evolution of the system. Indeed, the former are analytic functions, whereas the latter is not square integrable. The additional presence of a continuous spectrum in $P(x,D)$ is fundamental to explain this discrepancy.

\section{Conclusions}
\label{sec:conclusions}

In this work, we studied the relationship between internal wave attractors and the spectra of a class of zeroth-order pseudo-differential operators.

First, we developed numerical techniques to approximate the solution to the nonlocal wave equation \eqref{eq:evol_problem} and to the elliptic eigenvalue problem \eqref{eq:ev_problem}. The resulting methods are fourth order accurate in time and (locally) spectrally accurate in space. Given that the solution to the evolution problem develops attractors (singularities), global spectral accuracy cannot be expected. Moreover, RED estimates confirm that, as $t \to \infty$, the system evolves into a state that is not square-integrable. 

Then, we used these methods to analyze further spectral properties of the pseudo-differential operators in study. We showed that by appropriately reordering the viscous eigenvalues, we can capture $C^\infty$ trajectories that approximate the embedded eigenvalues as the viscosity decreases. Also, we compared side-by-side low-viscosity eigenmodes and the evolution (both in real and frequency space). Here, we related the observed behaviour to the geometrical structure of the energy surfaces on which the flows take place. In conclusion, the embedded eigenmodes do describe (at least, partially) the long-term dynamics of the problem.

Nevertheless, the viscous approximation of embedded eigenvalues and eigenfunctions still constitutes a challenging problem, given that as the viscosity decreases, so does the regularity of the viscous eigenfunctions. While taking more wave numbers (i.e. a finer spatial mesh) might seem like a good idea, the fact that the eigenvalues tend to cluster as the viscosity decreases creates a conflicting situation. Further research in this area is needed.

\section*{Acknowledgments}

We would like to thank Maciej Zworski for bringing us this very interesting problem and for the helpful discussions in this regard.

Javier A. Almonacid thanks the financial support of Simon Fraser University through the Graduate Dean's Entrance Scholarship. Nilima Nigam thanks the support of the Natural Sciences and Engineering Research Council of Canada (NSERC).

\small
\bibliographystyle{siam}
\bibliography{references}

\end{document}